\documentstyle[amssymb]{amsart}
\newcommand{\stopproof}{\hfill \nobreak\medskip $\blacksquare$ \\
\hspace*{\fill}}
\newcommand{\dom}{\mbox{\rm dom}}

\newcommand{\AND}{\mbox{ \rm and }}

\newcommand{\forces}[2]{\Vdash_{#1} \mbox{``} #2 \mbox{''}}

\newcommand{\proof}{{\bf Proof:} \ }

\newcommand{\PP}{{\Bbb P}}

\newcommand{\RR}{{\Bbb R}}

\newcommand{\presup}[2]{\, ^{#1} \! #2}

\newcommand{\fomom}{\presup{\omega}{\omega}}

\newtheorem{theor}{Theorem}[section]

\newtheorem{defin}{Definition}[section]
\newtheorem{corol}{Corollary}[section]
\newtheorem{lemma}{Lemma}[section]
\newtheorem{quest}{Question}[section]

\author{Juris Stepr\={a}ns}
\address{Department of Mathematics, York University \\
4700 Keele Street \\
North York, Ontario \\ Canada \ \ \ \ M3J 1P3}
\newcommand{\success}{\mbox{succ}} 
\newcommand{\rank}{\rho} 
\newcommand{\poset}{\PP} 
\newcommand{\interior}{\mbox{int}} 
\newcommand{\concat}{\wedge} 
\newcommand{\sakne}{\mbox{root}} 
\newcommand{\BH}{\mbox{branching-height}} 
\newcommand{\rabove}[1]{\langle #1\rangle} 

\title{Cardinal Invariants Associated with Hausdorff Capacities}
\begin{document} 
\bibliographystyle{amsplain}
\maketitle

\begin{abstract}
Let $\lambda(X)$ denote Lebesgue measure.
If $X\subseteq [0,1]$ and $r \in (0,1)$ then the $r$-Hausdorff capacity of
$X$ is denoted by $H^r(X)$ and is defined to be the infimum of all $\sum_{i=0}^\infty
\lambda(I_i)^r$ where $\{I_i\}_{i\in\omega}$ is a cover of $X$ by
intervals. The $r$ Hausdorff capacity has the same null sets as the
$r$-Hausdorff measure which is familiar from the theory of fractal
dimension. It is shown that, given $r < 1$,  it is possible to enlarge a model of set theory, $V$,
by a generic extension $V[G]$ so that the reals of $V$ have Lebesgue measure zero
but still have positive $r$-Hausdorff capacity.
\end{abstract}

\section{Introduction} If $r \in [0,1]$ then for any set $X\subseteq [0,1]$ the $r$-Hausdorff
capacity of $X$ is denoted by $H^r(X)$ and is defined to be the
infimum of all $t$ such that there is a cover of $X$ by intervals, $X
\subseteq \bigcup_{i = 0}^\infty I_i$, such that $t = \sum_{i = 0}^\infty
\lambda(I_i)^r$. This notion may be familiar from its use along the way to defining $r$-Hausorff
measure. Given $\beta > 0$,  $H^r_\beta(X)$ is defined, for any set $X\subseteq [0,1]$,
to be the
infimum of all $t$ such that there is a cover of $X$ by intervals, $X
\subseteq \bigcup_{i = 0}^\infty I_i$, such that $t = \sum_{i = 0}^\infty
\lambda(I_i)^r$ and such that the length of each interval $I_i$ is less than $\beta$. The
$r$-Hausdorff measure of a set $X$ is then defined to be the supremum of $H^r_\beta(X)$ as $\beta$
ranges over all positive real numbers. However, the topic of this paper if $r$-Hausdorff capacity
rather than $r$-Hausdorff measure. The crucial difference between the two is that, while the
$r$-Hausdorff measure is countably additive, the $r$-Hausdorff capacity is only subadditive if $r
\in (0,1)$. 
A proof of the fact that $H^r$ is actually a capacity can
be found in \cite{rogers.1} on page 90. For more details on Hausdorff measure
consult \cite{rogers.1}, \cite{carleson} or \cite{della}.

For the rest of this paper let $r$ be a fixed real number such that $ 0 < r < 1$.
Let $\lambda(X)$ denote the Lebesgue measure of any measurable set $X\subseteq [0,1]^n$.
 It will be shown  that it is possible to genericall y extend  an arbitrary model of set theory
so that
the ground model reals have Lebesgue measure 0 but still have
positive $r$-Hausdorff measure.  If this process could be iterated $\omega_2$ times and the ground
model satisfied the Continuum Hypothesis then it would yield
a model every set of size $\aleph_1$ has Lebesgue measure zero yet there is a set of size $\aleph_1$
which has positive $r$-Hausdorff capacity. This raises the following conjecture which uses the
obvious extension of Hausdorff capacity to $\RR^n$  : For any $n \in \omega$ and $r < n$ it is
consistent
that every set of size $\aleph_1$ has Lebesgue measure zero yet there is a set $X\subseteq [0,1]^n$
of size $\aleph_1$ which has positive $r$-Hausdorff capacity.

This is related to the following
 question posed by P. Komjath.
 \begin{quest} \label{Komjath} Suppose that every set of size $\aleph_1$ has Lebesgue measure zero.
Does
it follow that the union of any set of $\aleph_1$ lines in the plane   has Lebesgue measure zero?
\end{quest}
To see the relationship between this question and $r$-Hausdorf capacity consider that it is easy to
find countably many unit sqaures in the plane such that each line passes through either the top and
bottom or the left and right sides of at least one of these squares. It is therefore possible to
focus attention on all lines which pass through the top and bottom of the unit square. For any such
line $L$ there is a pair $(a,b)$ such that both the points $(a,0) $ and $(1,b)$ belong to $L$. If
the mapping which send a line $L$ to this pair $(a,b)$ is denoted by $\beta$ then it is easy to see
that $\beta$ is continuous and that if $S\subseteq [0,1]^2$ is a square of side $\epsilon$ then the
union of $\beta^{-1}S$ has measure $\epsilon$ while $S$ itself has measure $\epsilon^2$. In other
words, the Lebesgue measure of the union of $\beta^{-1}X$ is no larger than the 1-Hausdorff capacity
of $X$ for any $X\subseteq [0,1]^2$. Hence, if the answer to Question~\ref{Komjath} was negative
this would imply that  the conjecture is true. While this was the motivation for studying    the
problem of Hausdorff capacity, it may  be that the notion of Haudsorff capacity is actually more
central than Komjath's question itself.

If $1\leq n\leq m$ then
define $\pi : [0,1]^m \to [0,1]^n$ by $$\pi_n(x_1, x_2, \ldots,x_m) =
(x_1, x_2, \ldots,x_n).$$ If  $X \subseteq [0,1]^m$ then $\pi_n(X)$
will denote the image of $X$ under the mapping $\pi_n$.                     If $A\subseteq
[0,1]^{d}$ and $1 \leq n < d$ then, for any $x\in [0,1]^n$, the notation $A_x$ will be used to
denote $\{y \in [0,1]^{d-n} : \pi_n(x,y) = x\AND (x,y) \in A\}$.

\section{A General Class of Forcing Partial Orders}
This section will be devoted to examining a generalization of
Superperfect forcing obtained by insisting that on a dense set of
nodes the splitting is into a set of positive measure with respect to
some ideal. Such generalizations have been considered by various
authors.   Throughout this section the term ideal will always refer to a proper ideal on $\omega$
which contains all finite sets. In later sections ideals will be constructed on countable sets other
than $\omega$, but it will simplify notation to ignore this for now.

\begin{defin}
Let ${\cal I} = \{{\cal I}_n\}_{n \in \omega}$ be a sequence of ideals.
  The partial
order $\poset({\cal I})$ will be defined to consist of trees
$T\subseteq \bigcup_{n\in\omega} \prod_{i\in n}D_i$ such that for
every $t \in T$ one of the following two alternatives holds:
\begin{itemize}
\item $|\{n\in \omega : t\concat n \in T\}| = 1$
\item $\{n\in \omega : t\concat n \in T\} \in {\cal I}_{|t|}^+$

\end{itemize} If $\{n\in \omega : t\concat n \in T\} \in {\cal
I}_{|t|}^+$ then $t$ will be said to be a branching node of $T$ and
the set of branching nodes will be denoted by $B(T)$. Define
$\poset({\cal I})$ to consist of all $T$ such that for every $t \in T$
there is $s \in B(T)$ such that $ t \subseteq s$.  The ordering on
$\poset({\cal I})$ is inclusion.
\end{defin}
 It is  left to the reader to verify that
$\poset({\cal I})$        is                                                      proper and,
indeed,  that  it satisfies Axiom A. A standard argument works.

Suppose now that $T \in \poset({\cal I})$. Then the root of $T$ is the
unique minimal member of $B(T)$ and is denoted by $\sakne(T)$.  If $t
\in B(T)$ then the set of successors of $t$ is denoted by
$\success_T(t)$ and is defined by $\success_T(t) = \{n\in \omega :
t\concat n \in T\}$.  The branching height of $t$ will be denoted by
$\BH(t)$ and is defined to be $|\{s\subseteq t : s\in B(T)\}|$ --- so
$\BH(\sakne(T)) = 1$.  A subset $S\subseteq T$ will be said to be a
{\em subtree} if it is closed under taking initial segments. The tree
{\em generated} by $X \subseteq T$ is simply the set of all initial segments
of members of $X$. Observe that $\success_S$ can be defined for any
subtree, regardless of whether or not $S\in\poset({\cal I})$.  A
subset $S\subseteq T$ will be said to be a {\em full subtree} of $T$
if and only if for every $t \in S$ either $t$ is a maximal member of
$S$ or $\success_T(t) = \success_W(t)$.  If $t \in T$ then
$T\rabove{t}$ is defined to be the subtree of $T$ consisting of all $s
\in T$ such that either $s\subseteq t$ or $ t \subseteq s$. If $S
\subseteq T$ is a  subtree then define the {\em interior}  of $S$ to be
the set of all non maximal elements of $B(T) \cap S$ and denote this
by $\interior(S)$ --- the dependence on $T$ will suppressed.

If $T \in \poset ({\cal I})$ then define a function $\Psi$ on $B(T)$
to be {\em approximating} if $\Psi(t)\subseteq [0,1]$ is a finite
union of rational intervals for each $t \in B(t)$ and it is monotone
in the sense that $\Psi(t) \subseteq \Psi(s)$ if $t \subseteq s$. If
$T \in \poset ({\cal I})$, $ x\in [0,1]$ and $\Psi$ on $B(T)$ is
approximating then define $R(T,\Psi, x)$ to be the tree generated by
$\{t \in B(T) : x \notin \Psi(t)\}$.

\begin{defin} An ideal \label{KeyProperty} ${\cal I}$ will be said to satisfy KP$(r)$ if and only
if for all
\begin{itemize} \item $\theta  < 1$
\item
$X \in {\cal I}^+$ \item functions $F$ from $X$ to the Borel  subsets of
$[0,1]$ satisfying that $H^r(F(x))
\leq \theta$  for each $x\in X$ \item $\epsilon > 0$ \end{itemize} there is
some $Y\subseteq [0,1]$ as well as $Z\subseteq [0,1] $ such that
\begin{itemize} \item $H^r(Y) \leq
\theta$ \item $\lambda(Z) < \epsilon$ \item $\{x \in X :  y \notin F(x)\} \in
{\cal I}^+$ for every $y \in [0,1] \setminus (Y\cup Z)$
\end{itemize}\end{defin}

Well founded trees will play an important role in the following
discussion but the standard equivalance between well founded trees and
trees with rank functions is not as convenient a slight modification
of this notion.
If $T \in \poset({\cal I})$ and $S\subseteq T$ then
the {\em standard rank} of $S$ will denote the rank of $S\cap B(T)$
when this is considered as a tree under the inherited
ordering. Later on, a different rank function will be introduced and it should not be confused with
the standard rank.

If $T\in \poset({\cal I})$ and $W\subseteq T$ is a full subtree then $W'\subseteq W$ will be said to
be {\em large} if:
\begin{itemize}
\item $\sakne(W) \in W'$
\item if $t \in W'\setminus B(T)$ then $t$ is not maximal in  $W'$
\item if $t \interior(W)\cap W'$ then $\success_{W'}(t) \in {\cal I}_{|t|}^+$.
\end{itemize}

\begin{lemma}
Suppose that
\begin{itemize}
\item  ${\cal I} = \{{\cal I}_n : n\in \omega\}$ is a sequence of ideals,
each satisfying KP$(r)$
\item $T\in \poset({\cal I})$\label{man.2} and
$W\subseteq T$ is a well founded full subtree of standard rank $\beta$
\item  $\Psi$
is an approximating function on $B(T)\cap W $
\item $\theta < 1$
\item
$H^r(\Psi(t) )<\theta$ for any $t
\in  W\cap B(T)$\end{itemize}
then, there is some $x\in [0,1]$ such that $R(W,\Psi, x)$ is a large subtree
of $W$.
\end{lemma}
\proof   It will be shown by induction on $\beta \in \omega_1$ that the following, stronger
condition holds: \medskip

\noindent {\bf Q}($\beta$): \ \ \  If  $ s \in T$ and $W\subseteq T\rabove{s}$ is a well founded
full subtree
of standard rank $\beta$, $\theta < 1$, $\epsilon > 0$ and $\Psi$ is an approximating function
on $B(T)\cap W $ such that $H^r(\Psi(t) )<\theta$ for any $t\in W\cap
B(T)$ then, there is some $Y\subseteq [0,1]$ and $Z\subseteq [0,1]$
such that
\begin{itemize}
\item for each $x \in [0,1] \setminus ( Y\cup Z)$
 $R(W,\Psi,x)$ is a large subtree of $W$
 belongs to ${\cal I}^+_{|t|}$
\item $H^R(Y) \leq \theta$
\item $\lambda(Z) <\epsilon$.
\end{itemize}

First notice that this implies the lemma by choosing  $s = \sakne(T)$ and
$\epsilon < 1 -
\theta$ because $\lambda(Y) \leq H^r(Y)$ and so, if $\lambda(Z) <
\epsilon$ then $[0,1] \setminus (Y\cup Z)$ is not empty.
If $\beta = 0$ the statement is vacuous and
  Q(1) is implied by KP($r$). Now
assume that Q($\gamma$) has been established for all $\gamma \in
\beta$. If $W$ is a well founded
subtree of $T\rabove{s}$, $\theta < 1$, $\epsilon > 0$ and $\Psi$ is an
approximating function on $W\cap B(T)$ then for each $ n \in
\success_T(s)$ the standard rank of $W\rabove{s\concat n}$ is less than $\beta$.
Moreover, $H^r( \Psi(t) )<\theta$ for any $t \in W\rabove{s \concat n}\cap
B(T)$.  It is therefore possible to apply the induction hypothesis to $T\rabove{s\concat n}$ for
each $n\in \success_T(s)$ to find $Y_n$ and $Z_n$ such that
\begin{itemize}
\item $H^r(Y_n) \leq \theta/2^{n+1}$
\item $\lambda(Z_n) < \epsilon/2^{n+1}$
\item if $ x\in [0,1]\setminus (Y_n \cup Z_n )$
 then $R(W\rabove{s\concat n}, \Psi, x)$ is a large subtree of $W\rabove{s\concat n}$.
\end{itemize}

Now let  $X =  \success_T(s) \in {\cal
I}_{|s|}^+$.
Choose  a function $F$
on $X$ such that  $F(d) \supseteq Y_{d}$ and $F(d)$ is a $G_\delta$
such  that $H^r(F(d)) = H^r(Y_{d})
\leq \theta$ for each $ d \in X$.  It follows from KP($r$) that  there are $Y'\subseteq [0,1]$ and
$Z'\subseteq [0,1]$ such that \begin{itemize}
\item $H^r(Y') \leq \theta/2$
\item $\lambda(Z') < \epsilon/2$
\item $\{d \in X : x \notin F(d)\} \in {\cal I}_{|s|}^+$ for each $ x\in [0,1]\setminus (Y \cup
Z')$ \end{itemize} Now let $Z= Z' \cup \cup\{Z_n : n \in
X \}$ and note that $\lambda(Z) <
\epsilon$ and, by the subadditivity of $H^r$, $H^r(Y) \leq \theta$. Hence, in order to verify that
Q($\beta$)
holds it suffices to show that if $x \in [0,1] \setminus (Y\cup Z)$ and $ t \in R(W,\Psi,x)\cap
\interior(W)$ then $\success_{R(W, \Psi,x)}(t) \in {\cal I}^+_{|t|}$. If $t = s$
this follows from the application of KP($r$) and the fact that $x\notin Y'\cup Z'$. In every
other case it follows from the use of the induction hypothesis because $t \supseteq s\concat n$ for
some $n$ and, therefore, $ x \notin Y_n\cup Z_n$ implies that   $\success_{R(W, \Psi,x)}(t) \in {\cal I}^+_{|t|}$.
\stopproof

For the remainder of this section fix a sequence  of ideals ${\cal I} = \{{\cal I}_n : n\in
\omega\}$ and $T
\in \poset({\cal I})$. For $t \in B(T)$ and $n \in \success_T(t)$ define $t\oplus   n$ to be the least $s
\in B(T)$ such that $t\concat n \subseteq s$.  If $X\subseteq
T$ is a subtree then a rank function $\rank_{X}$ will be defined on $B(T)\cap X$ by bar
 induction.
To begin, define $\rank_{X}(t)$  to be 0 if
there is some $t'\subseteq t$ such that $t' \in B(T)$ and $\success_X(t')
 \in {\cal I}_{|t'|}$.
Define $\rank_{X}(t) $ to be the least ordinal $\beta$ such that there is some
$A \in {\cal I}_{|t|}$    such that $\rank_{X}(t \oplus    n)$ is defined for each $ n \in
\success_X(t)\setminus A$ and $\rank_{X}(t \oplus    n) \in \beta$ for any such $n$.
The rank of $X$ is defined to be the rank of  its root, provided this is defined.
A subtree $X\subseteq T$  will be defined to be {\em small} if
 $\rank_{X}(t)$ is defined for all $ t \in B(T)\cap X$ and $\rank_X(t) > 0$ if and only if $ t \in
\interior(X)$.

\begin{lemma}
If
 $X\subseteq T$ is a subtree and there is some $t\in X\cap B(T)$ for which $\rank_X(t)$ is not
defined \label{manu.rank} then $X$ contains a member of $\poset({\cal I})$.
\end{lemma}
\proof  This is standard. Let $S$ be the subtree of $T$ generated by all
 $t \in X\cap  B(T)$ such that $\rank_{X}(t)$ is not defined.
Notice that if $t \in S$ then $$ \{n \in \success_X(t) :
\rank_{X}(t \oplus    n) \mbox{ is not defined}\} =
\success_X(t) $$
and this  belongs to $ {\cal I}_{|t|}^+$.  Hence $S \in \poset({\cal I})$
provided that it is not empty.  The hypothesis of the lemma guarantees
that this is not the case.
\stopproof

For any subtree $W  \subseteq T$ and any function $\theta \in \prod_{
w\in \interior(W)}{\cal I}_{|w|}$ define $$W^\theta = \{w \in W :(\forall n\in
\omega)(w(n) \notin
\theta(w\restriction n))\}$$ or, in other words,  $W^\theta$ is obtained by throwing away
${\cal I}_{|t|}$ many successors, determined by $\theta$,  of each  $t \in \interior(W)$.
If $W$ and $X$ are subtrees of $T$ define $W \prec X$ if and only if there exists $\theta \in
\prod_{w\in \interior(W)}{\cal I}_{|w|}$ such  that for every
$\theta' \in \prod_{x\in \interior(X)}{\cal I}_{|x|}$ there is a one-to-one function $G:W^\theta\cap
B(T) \to (X^{\theta'}\cap
B(T))\setminus \{\sakne(X)\}$ which is order preserving in the sense that if $ t\subseteq s$ then
$G(t) \subseteq G(s)$.

For any small subtree $Y \subseteq T$ a function $\theta \in \prod_{y\in \interior(Y)}{\cal
I}_{|y|}$ will be said
to be {\em  a witness to the rank} of $Y$ if and only if for each    $ t \in \interior(Y)$
$$\rank_{Y,T}(t\oplus   n) \in \rank_{Y,T}(t)$$ for each $ n \in \success_Y(t)\setminus
\theta(t)$.

\begin{lemma}
Let  \label{ctsfcn.8}
  $W$ and $ X$ be small subtrees of $T$
of rank $\alpha$ and $\beta$ respectively. If $ \alpha \in \beta$, then for any $\theta \in
\prod_{w\in \interior(W)}{\cal I}_{|w|}$ which is a witness to the rank of $W$ and any
 $\theta \in
\prod_{x\in \interior(X)}{\cal I}_{|x|}$ there is a one-to-one function $G_{\theta,\theta'}
:W^\theta\cap B(T) \to X^{\theta'}\cap B(T)$ which is order preserving such that
$G_{\theta,\theta'}(\sakne(W)) \neq \sakne(X)$.
Moreover, $G_{\theta,\theta'}$ is continuous in the variable $\theta'$ in the sense that
the mapping $\theta' \mapsto G_{\theta,\theta'}$ is a  continuous function from $
\prod_{x\in \interior(X)}{\cal I}_{|x|}$ to ${}^{W^{\theta}\cap B(T)}(X^{\theta'}\cap B(T))$
where ${\cal I}_{n}$ is considered as a subspace of $2^{\omega}$. \end{lemma}
\proof
Suppose that  $\alpha \in \beta$ and $W$ and $X$ are small subtrees
of $T$
of rank $\alpha$ and $\beta$ respectively.  Let $\theta \in
\prod_{w
\in \interior(W)}{\cal I}_{|w|}$ be a witness to the rank of $W$. For every $\theta' \in \prod_{
x\in \interior(X)}{\cal I}_{|x|}$ a function
$G_{\theta,\theta'}
:W^\theta\cap B(T) \to X^{\theta'}\cap B(T)$ can be defined by induction on the branching height of
nodes of
$W^\theta\cap B(T)$.
The induction hypothesis will be that $\rho_X(G_{\theta,\theta'}(t)) \geq \rho_W(t)$.
Define $G_{\theta,\theta'}(\sakne(W)) = \sakne(X)\oplus   m$ where $m$ is
the least integer such
that         $m \in \success_X(\sakne(X))\setminus \theta'(\sakne(X))$ and
$\rank_{X}(\sakne(X)\oplus   m) \geq \alpha$. Such an $m$ must exist because
$\rank_{X}(\sakne(X)) = \beta> \alpha$ and $\theta'(\sakne(X)) \in {\cal
I}_{|\sakne(X)|}$. If $t$ and $t\oplus   n$
are
both in $W^\theta\cap B(T)$ and $G_{\theta,\theta'}(t)$ and $G_{\theta,\theta'}(t\oplus   i)$ are defined
for
$ i \in n$ then define  $G_{\theta,\theta'}(t\oplus   n) = G_{\theta,\theta'}(t)\oplus   k$ where $k$ is the least
integer such that                               $$k\in \success_X(G_{\theta,\theta'}(t))\setminus
(\theta'(G_{\theta,\theta'}(t))\cup \{G_{\theta,\theta'}(t\oplus   i)(|G_{\theta,\theta'}(t)|) : i \in
n\})$$ and $\rank_{X}(G_{\theta,\theta'}(t\oplus   k)) \geq
\rank_{W}(t\oplus   n)$. The reason such a $k$ exists is that, by the  induction hypothesis,
$\rank_{X}(G_{\theta,\theta'}(t)) \geq
\rank_{W}(t)$ and hence $\rank_{X}(G_{\theta,\theta'}(t) >
\rank_{W}(t\oplus   n)$ because $\theta$ is a witness to the rank of $W$ and $t\oplus   n \in
W^\theta$. Since ${\cal I}_{|G_{\theta,\theta'}(t)|}$ contains all finite subsets it must be
that there is some   $$k\in \success_X(G_{\theta,\theta'}(t))\setminus
(\theta'(G_{\theta,\theta'}(t)) \cup \{G_{\theta,\theta'}(t\oplus   i)(|G_{\theta,\theta'}(t)|) : i \in
n\})$$ such that $\rank_{X}(G_{\theta,\theta'}(t\oplus   k) \geq
\rank_{W}(t)\oplus   n) $. Since the inductive
hypothesis is preserved, this construction can be carried out for all nodes in $W^\theta$.
Obviously $G_{\theta,\theta'}$ is a one-to-one, order preserving function. By construction,
$G_{\theta,\theta'}(\sakne(W)) \neq \sakne(X)$.

The continuity of $G_{\theta,\theta'}$ in  the variable $\theta'$ follows from the minimal choice of
the integer $k$
 such that                               $$k\in \success_X(G_{\theta,\theta'}(t))\setminus
(\theta'(G_{\theta,\theta'}(t))\cup \{G_{\theta,\theta'}(t\oplus   i)(|G_{\theta,\theta'}(t)|) : i \in
n\})$$ and $\rank_{X}(G_{\theta,\theta'}(t\oplus   k) \geq
\rank_{W}(t\oplus   n)$.  An open neighbourhood of $G_{\theta,\theta'}$ in
$${}^{W^\theta\cap B(T)}(X^{\theta'}\cap B(T))$$ is
specified by a restriction of $G_{\theta,\theta'}$
to a finite subset. Given a finite subset $a \subseteq W^{\theta}$ it is possible to
find a finite set $b \subseteq X^{\theta'}$ such that if $t\in a$, $G_{\theta,\theta'}(t)
= s\oplus   m$ and $i \in m \setminus \theta'(s)$ then  $s\oplus   i$ also belongs to $b$. Let $M\in
\omega$ be such that the range of each $t\in b$ is contained in $M$. It is easy to check that if
$\theta'' \in \prod_{x\in \interior(X)}{\cal I}_{|x|}$ is such that $\theta''(t)\cap M =
\theta'(t)\cap M$ for each $ t\in b$ then $G_{\theta,\theta''}\restriction a =
 G_{\theta,\theta'}\restriction a$.
\stopproof

          \begin{lemma}
If $W$ and $X$ are
small subtrees \label{whew} of $T$
of rank $\alpha$ and $\beta$ respectively, then $W \prec X$ if and only if $\alpha \in \beta$.
\end{lemma}
One direction is an immediate consequence of Lemma~\ref{ctsfcn.8} because a witness to the rank of
$W$ can always be found.
For the other, it will be shown by induction on
$\alpha$ that if
$\alpha \leq \beta$  then $X \not\prec W$. For $\alpha = 0$ this is trivial so assume that the
assertion has been established for all $\alpha' \in \alpha$ and suppose that
$X \prec W$. This means that there is some $\theta \in \prod_{x\in \interior(X^*)}{\cal I}_{|x|}$
such that
 for every $\theta' \in \prod_{
w\in \interior(W)}{\cal I}_{|w|}$ there is a one-to-one function
$G: X^\theta\cap B(T) \to W^{\theta'}\cap
B(T)$ which is order preserving such that $G(\sakne(X)) \neq \sakne(W)$.
Let $\theta'$ be a witness to the rank of $W$ and let the function
$G$ from $X^\theta\cap B(T)$ to $W^{\theta'}\cap
B(T)\setminus \{\sakne(W)\}$ be one-to-one and order preserving.
It must be  that $\rank_{W}(G(\sakne(X))) \in \rank_{W}(\sakne(W)) =
\alpha\leq \beta$. Therefore, it is possible to find $m\in \success_X(\sakne(X)) \setminus
\theta(\sakne(X))$ such that $\rank_{X}(\sakne(X) \oplus   m) \geq
   \rank_{W}(G(\sakne(X))) $. Obviously $G$, $\theta\restriction X\rabove{\sakne \oplus    m}$ and
$\theta'\restriction W\rabove{G(\sakne(X)}$ establish that $X\rabove{\sakne \oplus    m}    \prec
W\rabove{G(\sakne(X)}  $ which contradicts the induction hypothesis.
\stopproof

An ideal ${\cal I}$  will be said to be $\Sigma^1_1$ if it is a
$\Sigma^1_1$ subset of $2^\omega$ under the natural identification. The next lemma shows that the
relation $\prec$ is $\Sigma^1_1$ provided
that each of the ideals of ${\cal I}  = \{{\cal I}_n\}_{n\in\omega}$ is $\Sigma^1_1$. This will
require the full conclusion of Lemma~\ref{ctsfcn.8} since the obvious calculation only shows that
$\prec $ is $\Sigma^1_3$.

\begin{lemma}
If       each of the ideals of ${\cal I}  = \{{\cal I}_n\}_{n\in\omega}$ is $\Sigma^1_1$
        \label{sigma32sigma1} then the relation $\prec$ defined from them is $\Sigma^1_1$.
                              \end{lemma}
\proof    Since each ${\cal I}_n$ is $\Sigma^1_1$ it is possible to choose continuous functions
$f_n :\fomom\to 2^{\omega}$ such that ${\cal I}_n $ is the image of $f_n$.
First it will be shown that $W \prec X$ if and only if there is some $\theta\in \prod_{w\in
\interior(W)}{\cal I}_{|w|}$ and a continuous function
$$G :  \prod_{x\in
\interior(X)}\fomom \to {}^{W^\theta\cap B(T)}\theta\cap B(T)(B(T)\cap X)$$
such that \begin{enumerate}
\item for all $\theta' \in \prod_{x\in
\interior(X)}\fomom$, if $t$  belongs to $W^\theta\cap B(T)$ then
$G(t) = s \oplus  m$ for some $s \in B(T)cap  X$ and some  $m\in
\success_X(s)(t)\setminus f_{|s|}(\theta'(s))$
\item $G(\theta')$ is order preserving                               for all $\theta' \in
\prod_{x\in \interior(X)}\fomom$
\item $G(\theta')$ is one-to-one    for all $\theta' \in
\prod_{x\in
\interior(X)}\fomom$
\item $G(\theta')(\sakne(W)) \neq \sakne(X)$ for all $\theta' \in \prod_{x\in
\interior(X)}\fomom$
\end{enumerate}
Assuming $W \prec X$, it is possible to use Lemma~\ref{ctsfcn.8} to define $G(\theta') =
G_{\theta,\mu(\theta')}$ where $\mu(\theta')(x)
= f_{|x|}(\theta'(x))$ for each $ x\in \interior(X)$. Note that $G$ is continuous because of the
final sentence
of Lemma~\ref{ctsfcn.8} and the continuity of $\mu$, which is a consequence of the continuity of
each
$f_n$.     The other direction of the equivalence is clear because each $f_n$ is onto ${\cal I}_n$.

Hence it suffices to check that the clauses (1) - (4) are arithmetic.
Since $T$ and $B(T)$ can be used as parameters, the only problematic part
is the use of the quantifiers $$\mbox{for all } \theta' \in
\prod_{x\in
\interior(X)}\fomom \  .$$ However, the continuity of $G$ allows these to be replaced by quantifiers
over approximations to $\theta'$. In particular, it suffices to choose a countable dense subset
 $ C \subseteq \prod_{x\in
\interior(X)}\fomom$    and replace each instance of
$$\mbox{for all } \theta' \in
\prod_{x\in
\interior(X)}\fomom$$
         by  ``for all $\theta' \in C'$''.
\stopproof

\begin{lemma} For all $\alpha$ such    \label{largess}
that $1 \leq\alpha \in \omega_1$ and $t\in B(T)$ there is a
well founded full subtree $W$ of standard rank $\alpha$ such that $\sakne(W) = t$ and if
$W'\subseteq W$ is any  large subtree
 then $\rank_{W'}(t) = \alpha$.\end{lemma}
\proof Proceed by induction on $\alpha$. The case $\alpha = 1$ is trivial so assume the assertion
has been established for all $\alpha' \in \alpha$.
         First suppose that $\alpha = \beta +1$.    Given $t\in B(T)$ choose for each $n\in
\success_W(t)$ a well founded full subtree   $W_n$ of standard rank $\beta$ such that $\sakne(W_n) =
t\oplus   n$ and if $W'\subseteq W_n$ is any  large subtree
then $\rank_{W'}(t\oplus   n) = \beta$. Let $W = \cup_{n\in \success_W(t) }W_n$.

If $\alpha$ is a limit let $\{\beta_n\}_{n\in\omega}$ converge to $\alpha$ from below.
Given $t\in B(T)$ choose for each $n\in
\success_W(t)$ a well founded full subtree   $W_n$ of standard rank $\beta_n$ such that $\sakne(W_n)
= t\oplus   n$ and if $W'\subseteq W_n$ is any  large subtree
then $\rank_{W'}(t\oplus   n) = \beta_n$. Let $W = \cup_{n\in \success_W(t) }W_n$.
Since ${\cal I}_{|t|}$ contains all finite sets this works.
                                                   \stopproof

\section{The preservation Theorem}
\begin{theor} If ${\cal I} = \{{\cal I}_n : n\in \omega\}$ is\label{manu.main}  a sequence of
$\Sigma^1_1$ ideals satisfying KP$(r)$ such that each ${\cal I}_n $ contains all finite sets and
$G$ is $\poset({\cal I})$ generic over $V$ then $H^r([0,1]\cap V) = 1$ in $V[G]$.
\end{theor}
\proof  Suppose the theorem false --- in other words, there is some $\theta < 1$ and $\{J_n\}_{n =
0}^\infty$, a name for a sequence of intervals with rational
endpoints, as well as a condition $T\in \poset({\cal I})$ such that
$$T\forces{}{([0,1]\cap V) \subseteq \bigcup_{n=0}^\infty J_n}$$ and
$T\forces{}{\sum_{n=0}^\infty \lambda(J_n)^r < \theta}$.  By
thinning down $T$ it may be assumed that if $t \in B(T)$ and $\BH_T(t)
= n$ then $T\rabove{t}\forces{}{ J_i = J(t,i)}$ for some interval
$J(t,i) $ with rational endpoints, for each $ i\in n$.  Let $\Psi(t) =
\cup\{ J(t,i) : i \in |t|\}$ for each $t \in B(T)$.

 If there is some $x \in [0,1]$ such that $R(T,\Psi,x)$ contains some
$S \in \poset({\cal I})$ then it follows that $$S\forces{}{x\notin
\bigcup_{n=0}^\infty J_n}$$ contradicting the fact that $S \subseteq
T$. Hence by Lemma~\ref{manu.rank} it follows that
$\rank_{R(T,\Psi,x)}(t)$
is defined for all $ t \in R(T,\Psi,x)\cap B(T)$ and
therefore
$$T_x = R(T,\Psi,x)\setminus \{t \in R(T,\Psi,x) :(\exists
t'\subsetneq t)
\rank_{R(T,\Psi,x)}(t) = 0\}$$ is a small subtree for
each $ x\in [0,1]$. Notice that ``$\rank_{R(T,\Psi,x)}(t) =
0$''
is a $\Sigma^1_1$ statement because each ${\cal I}_n$ is $\Sigma^1_1$.
Hence
$\{T_x : x\in [0,1]\}$ is a $\Sigma^1_1$ set.

Since the relation $\prec$ defined on $\{T_x : x\in [0,1]\}$ is also $\Sigma^1_1$ by
Lemma~\ref{sigma32sigma1}, it follows from the Kunen-Martin Theorem
\cite{mosc}  and Lemma~\ref{whew} that there is some $\alpha \in
\omega_1$ such that the rank of $T_x$ is less than $\alpha$ for each
$x \in [0,1]$.  Use Lemma~\ref{largess}  to find $W \subseteq T$,  a well founded full subtree of
$T$,  of standard  rank $\alpha $
 such that  if
$W'\subseteq W$ is any large subtree
 then  the rank of $W'$ is $ \alpha $.
Observe that $\Psi$ is an approximating  function on
$B(T)\cap W $ such that $H^r( \Psi(t)
)<\theta$ for any $t
\in  W\cap B(T)$.
It follows from Lemma~\ref{man.2} that there is some $x\in [0,1]$ such
that  $R(W,\Psi,x)$ is a large subtree of $W$.
It follows that the rank of $R(W,\Psi,x)$ is at least $\alpha $ and this is
a contradiction because it implies that the rank of $T_x$ is at least
$\alpha$.
\stopproof

The reasonable conjecture at this point is that the conclusion of
Theorem~\ref{manu.main} holds for the countable support iteration of
the partial order $\poset({\cal I})$ . A proof of this would require
modifying the
preservation technology of Judah-Shelah cite{Juda.Shel, Ba.Ju.Sh} which was originally developed
to show that certain iterations preserve that the ground model reals
have positive Lebesgue measure.

\section{The Relation $\Xi$}
Sets with positive $r$-Hausdorff capacity may have measure zero but this type of set will play an
important role in the following discussion.  One would like to be able to say that  if
$\lambda(X) > 0$ then $H^r(X)$ can be calcultaed from $\lambda(X)$ or, at the very least, one
would might hope for  some relationship between $H^r(X)$ and  $\lambda(X)$. There are easy
counterexamples to this though. Let $X$ be such that  $H^r(X) = h > 0$ and $\lambda(X) = 0$ and then
$\lambda(X \cup [0,a]) = a$    and note that there is obviously no connection between $H^r(X\cup
[0,a])$ and $a$ when $a$ is much smaller than $h$. This sort of example is eliminated by
introducing a relation on sets which, roughly speaking, calculates the infimum of $H^r(X\setminus
Z)$ as $Z$ ranges over set of small measure.  The $\Xi$ relation, which is introduced in the next
definition, expands on this.

\begin{defin}
If $X$ and $Y$ are subsets of $[0,1]$ then define the relation
$\Xi_{\delta,\epsilon}(X,Y)$ to hold if and only if for every set $Z$, if
$\lambda(Z) < \epsilon$ then $H^r(X\cap Y \setminus Z) > H^r(Y) - \delta$.
If $X$ and $Y$ are subsets of $[0,1]^{n+1}$ then define the relation
$\Xi_{\delta,\epsilon}(X,Y)$ to hold if and only if
$$\Xi_{\delta,\epsilon}(\{x\in \pi_1(Y) :
\Xi_{\delta,\epsilon}(X_x,Y_x)\}, \pi_1(Y)).$$
\end{defin}
 The
relation $\Xi_{\delta,\epsilon}$ on $[0,1]^n$ can be considered as  a crude substitute for an
integral  when $n > 1$. In fact, one  might be tempted to define a better approximation to an
integral in the following way.
Define $\Xi'_\epsilon(A) = \inf\{H^r(A\setminus Z) : Z\subseteq [0,1] \AND \lambda(Z) < \epsilon\}$
for $A\subseteq [0,1]$. If $A\subseteq [0,1]^{n+1}$ then define $$\Xi'_\epsilon(A) = \sup\{\delta :
\Xi'_\epsilon(\{x\in [0,1] : \Xi'_\epsilon(A_x) \geq \delta\}) \geq \delta\}$$ by induction on $n$.
Notice however that the inequality $\Xi'_\epsilon(A) > \Xi'(B) - \delta$ is not equivalent to
$\Xi_{\delta,\epsilon}(A,B)$ even if $A\subseteq B$. The point is that
 if   $X$ and $Y$ are subsets of $[0,1]^{n+1}$ and $x \in \pi_1(Y)$ then it is possible
that $\Xi_{\delta,\epsilon}(X_x,Y_x)$ holds even though $ x\notin \pi(X)$.
This section collects
some facts about the $\Xi$ relation.

\begin{lemma}
If \label{ddimm} $A$ and $B$ are subsets of $[0,1]^n$ and $\Xi_{\delta,\epsilon}(A,B)$  holds then
for any $Z \subseteq [0,1]^n$ such that $$\lambda(Z) < \left (\frac{\epsilon}{2}\right )^n$$
$\Xi_{\delta,\epsilon/2}(A\setminus Z, B)$ also holds.\end{lemma}
\proof This is easliy proved using induction on $n$ and a simple
application of Fubini's Theorem.
\stopproof

The next two lemmas show that the relation $\Xi$ could have defined
from the top down rather than from the bottom up.
\begin{lemma}
If $A$ and $B$ are subsets of $[0,1]^{d+1}$ and \label{reverse}
$\Xi_{\delta,\epsilon}(\pi_d(A),\pi_d(B))$ holds and
$\Xi_{\delta,\epsilon}(A_x,B_x)$ holds for each $x \in \pi_d(A)$ then
$\Xi_{\delta,\epsilon}(A,B)$ holds as well.
\end{lemma}
\proof Proceed by induction on $d$ noting that the case $d=1$  follows
from the definition.
\stopproof

The proof of the next three lemmas are easy and left to the reader.

\begin{lemma} If $A$ and $B$ are subsets of $[0,1]^{d+1}$ and
\label{Xi.5} $\Xi_{\delta, \epsilon}(A,
B)$ holds then so does $\Xi_{\delta, \epsilon}(\{x\in \pi_d(B) : \Xi_{\delta,\epsilon}(A_x,
B_x)\}, \pi_d(B))\}$.
\end{lemma}

\begin{lemma} If $A$ and $B$ are closed subsets of $[0,1]^{d+1}$, $\delta > 0$ and
$\epsilon > 0$ then $$\{x \in \Pi_d(B) : \Xi_{\delta, \epsilon}(A_x,
B_x)\}$$ is a Borel
\label{Xi.2} set.
\end{lemma}

\begin{lemma} If
\label{Xi.3} $\Xi_{\delta, \epsilon}(A,
B)\}$ holds and $A\subseteq B$ then $\Xi_{\delta, \epsilon}(A,
A)\}$ holds as well.
\end{lemma}

\begin{lemma} If $\epsilon$, $\delta_1$ and $\delta_2$ are greater than 0,
\label{Xi.4}
$\Xi_{\delta_1, \epsilon}(A,
B)\}$  and $\Xi_{\delta_2, \epsilon}(B,
C)\}$ both hold and $B\subseteq C$ then $\Xi_{\delta_1 + \delta_2,
\epsilon}(A\cap
B), C\}$ also holds.
\end{lemma}

\begin{lemma}
If \label{ExTech} $D\subseteq [0,1]$ is a set such that
$\Xi_{\delta,\epsilon}(D,D)$ holds and $ \delta < H^r(D)$ then
there is some $\bar{\epsilon} > 0$
such that $\Xi_{\delta - \bar{\epsilon}, \epsilon/2}(D,D)$ holds as well.  \end{lemma}
\proof Let $\bar{\epsilon} < \min\{H^r(D) - \delta, r^{\frac{1}{1-r}},
\epsilon\}$.
Suppose that that $\lambda(Z) <
\epsilon/2$ and that $\{I_i :  i\in \omega\}$ is a cover of $D\setminus Z$.  Since
$\bar{\epsilon} < \theta$, by taking a tail of the sequence it is possible
to find $i_0\in \omega$ such that $\sum_{i = i_0 + 1}
^\infty\lambda(I_i) \leq \bar{\epsilon}/2  <
\sum_{i = i_0}^\infty\lambda(I_i) $ because, if $\sum_{i =
0}^\infty\lambda(I_i)\leq \bar{\epsilon}/2 < \epsilon/2$ then $\lambda(Z \cup
(\bigcup_{i\in\omega} I_i)) < \epsilon$ and so $Z \cup
(\bigcup_{i\in\omega} I_i)$ can not possibly be a cover of $D$.
  Let $J$ be an initial subinterval of $I_{i_0}$ such
that $\sum_{i = i_0+1}^\infty\lambda(I_i) + \lambda(J) = \bar{\epsilon}/2$ and
note that $\lambda(J) > 0$.  It follows that $\lambda(Z\cup
J\cup\bigcup_{i > i_0}I_i) <\epsilon$ and so $\lambda(I_{i_0}\setminus
J )^r +
\sum_{i \in i_0}\lambda(I_i)^r > \theta$ and hence, using the fact that $r < 1$,
$$\lambda(I_{i_0}\setminus J )^r +\lambda(J) + \sum_{i\neq
i_0}\lambda(I_i)^r \geq
\lambda(I_{i_0}\setminus J )^r +\lambda(J) + \sum_{i  \in
i_0}\lambda(I_i)^r + \sum_{i = i_0 + 1}^\infty\lambda(I_i) \geq
\theta + \bar{\epsilon}/2$$
Now notice that $$\sum_{i\in\omega} \lambda(I_i)^r =
\lambda(I_{i_0}\setminus J )^r +\lambda(J) + \sum_{i\neq i_0}\lambda(I_i)^r -
\left(\lambda(I_{i_0}\setminus J )^r + \lambda(J) -
\lambda(I_{i_0})^r\right ) $$
$$ \geq \theta + \bar{\epsilon}/2
-\left (\lambda(I_{i_0}\setminus J )^r + \lambda(J) -
\lambda(I_{i_0})^r\right) . $$ Hence
all that has to be shown is that $\lambda(I_{i_0}\setminus J )^r
+\lambda(J) - \lambda(I_{i_0})^r \leq
0$.

 To see this, define for $ a > 0$ $$F_a(x) = (a- x)^r + x
-a^r$$ and observe that $\frac{d}{dx}F_a(x) = 1 -
\frac{r}{(a-x)^{1-r}}$ and notice that this is negative if $ x < a <
\bar{\epsilon}$. Moreover, $F_a(0) = 0$ and so   $F_a$ is negative on the
interval $(0,a)$ if $a < \bar{\epsilon}$.
Because $0 <  \lambda(J) < \lambda(I_{i_0}) <
\bar{\epsilon}$ it follows that
$\lambda(I_{i_0}\setminus J )^r
+\lambda(J) - \lambda(I_{i_0})^r = F_{\lambda(I_{i_0})}(\lambda(J)) \leq 0$.
\stopproof

\begin{defin}
A subset $X\subseteq [0,1]$ will be said to be elementary if and only
if $$X = [p_0,q_0]\cup [p_1,q_1]\cup \ldots \cup[p_k,q_k] $$
where $p_i$ and $q_i$ are rational
numbers such that $p_i < q_i < p_{i+1}$ for each $i \in k$. A subset
$X \subseteq [0,1]^{n+1}$ is elementary if and only if there is an
elementary\label{elemntdef} set $ [p_0,q_0]\cup [p_1,q_1]\cup \ldots \cup[p_k,q_k]
\subseteq [0,1]$ such that $$X =  \bigcup_{i=0}^k [p_i,q_i]\times
X_i$$ and each $X_i\subseteq [0,1]^n$ is elementary.
\end{defin}

\begin{lemma}
If $U\subseteq [0,1]^n$ is open and $X\subseteq [0,1]^n$ is closed
then $\Xi_{\delta,\epsilon}(U,  X)$ if and only if for
\label{ReduceOpenToElementary} every $\bar{\epsilon} <
\epsilon$ there is an elementary
$Y\subseteq U$ such that $\Xi_{\delta,\bar{\epsilon}}(Y, X)$.
\end{lemma}
\proof    To begin, suppose that for every $\bar{\epsilon} < \epsilon$
 there is an elementary
$Y\subseteq U$ such that $\Xi_{\delta,\epsilon}(Y,X)$.
 Let
$Y_{k}\subseteq U$ be such that
$\Xi_{\delta,\frac{k\epsilon}{k+1}}(Y_{k}, X)$ and let $Y =
\cup_{ k>0} Y_{{k}}$. It will be shown by induction on $n$  that
$\Xi_{\delta,\epsilon}(U, X)$.

In particular, it will be shown by induction on $n$ that if $Y_m\subseteq [0,1]^n$ are sets such
that \begin{itemize}
\item each $Y_m$ is an elementary subset of $U$
\item $\Xi_{\delta, \epsilon(m)}(Y_m,X)$
\item $Y_m\subseteq Y_{m+1}$
\item  $\epsilon(m) \leq \epsilon(m+1)$
\item $\lim_{m\to \infty}\epsilon(m) = \epsilon$
\end{itemize}
then $\Xi_{\delta,\epsilon}(\cup_{m\in\omega}Y_m, X)$.
To begin the induction note that
the case $n=1$ is easy.                               Now suppose that the assertion has been
established for $n$ and that $X$ and $Y_m$ are subsets of $[0,1]^{n+1}$ for $m \in \omega$.
Suppose that  $Z\subseteq [0,1]$ is such that $\lambda(Z) < \epsilon$
and  $A = \{x\in \pi_1(X) : \Xi_{\delta, \epsilon(m)}(U_x,X_x)\}$
$$H^r(\{x\in \pi_1(X) : \Xi_{\delta, \epsilon}(A\setminus Z) > H^r(\pi_1(X)) - \delta$$
and define $$W_m = \{x\in \pi_1(X) : \Xi_{\delta, \epsilon(m)}((Y_m)_x,X_x)\}$$ for each $ m\in
\omega$. Note that if $\epsilon(m) > \lambda(Z)$ then $\lambda(W_n\setminus A) \geq \epsilon(m) -
\lambda(Z)$ because, otherwise, $\lambda(Z \cup (W_n\setminus A) ) < \epsilon(m)$ and so
$H^r(W_m \setminus ((W_m\setminus A)\cup Z) > H^r(X) - \delta$ contradicting that
                      $ W_m \setminus ((W_m\setminus A)\cup Z \subseteq A\setminus Z$.
Let $j$ be such that $\epsilon(j) > \lambda(Z)$. Then $\{W_m\setminus A : m > j\}$ is a a family of
measurable sets --- the measurablity follows from Lemma~\ref{Xi.2} and the fact that  each
$W_m$ is elementary, and hence, closed --- each of measure at least $\epsilon(j) - \lambda(Z) > 0$.
Hence there is some
$x\in [0,1]$ such that there are infinitely many $m \in \omega$ such that $x \in W_m \setminus A$.
Therefore there infinitely many $m \in \omega$ such that $\Xi_{\delta, \epsilon(m)}((Y_m)_x,X_x)$
and, by the induction hypothesis, it follows that $\Xi_{\delta, \epsilon}(U_x, X_x)$ holds because
$(Y_m)_x \subseteq (Y_{m+1})_x\subseteq U_x$, each $(Y_m)_x$ is elementary.

Conversely, suppose that  $\Xi_{\delta,\epsilon}(U,X)$ and let $\bar{\epsilon} <
\epsilon$.
Proceed by induction on $n$. The case $n=1$ is easy since $Y$ can be chosen to be a finite union of
intervals such that $\lambda(X\setminus Y) < \epsilon - \bar{\epsilon}$. Therefore suppose that
that    $U\subseteq [0,1]^{n+1}$  and $X\subseteq [0,1]^{n+1}$.  Using the induction hypothesis, it
follows that for each $x \in \pi(X)$,
if
 $
\Xi_{\delta,\epsilon}(U_x, X_x)$                 then
there is an elementary set $Y_x\subseteq U_x$ such that
$Y_x \subseteq_{\delta,\bar{\epsilon}}X_x$. Moreover, since $U$ is
open it follows that if $Y_x \neq \emptyset$ then there is
an open interval $J_x$ containing $x$ such that $J_x\times Y_x \subseteq U$. If $Y_x
= \emptyset$ then, since $X$ is closed there is an interval $J_x$
containing $x$
such that $\Xi_{\delta,\epsilon}(\emptyset, X_z)$ holds for all $z\in J_x$ .  Since
$\Xi_{\delta,\epsilon}(U, X)$ holds it must be the case that $$\Xi_{\delta,\epsilon}(\cup\{J_x :
\Xi_{\delta,\bar{\epsilon}}(Y_x, X_x)\}, \pi_1(X))$$
must hold. From the case $n=1$ it is possible to find an elementary set $Z   \subseteq
\cup\{J_x :
\Xi_{\delta,\epsilon}(U_x, X_x)\}$ such that        $\Xi_{\delta,\bar{\epsilon} + (\epsilon -
\bar{\epsilon})/2}(Z,\pi_1(X))$. Since $Z$ can be covered by finitely many intervals $J_x$
it is possible to obtain $Z' = \cup_{i\in j}J_i$ such that \begin{itemize}
\item $\lambda(Z\setminus Z') < (\epsilon - \bar{\epsilon})/2$
\item $J_i\cap J_{i'} = \emptyset$ if $ i \neq i'$
\item for each $i \in j$ there is some $x(i) $ such that $J_i \subseteq J_{x(i)}$
\end{itemize} It follows that $\Xi_{\delta,\bar{\epsilon}}(Z',\pi_1(X))$ so
let $Y = \cup_{i\in j} J_i\times Y_{x(i)}$.
     \stopproof

\begin{corol}
If $U$ is open \label{Xi.2.c} and $X$ is closed then the relation $\Xi_{\epsilon,\delta}(U,X)$ is
Borel. \end{corol}
\proof  It follows from Lemma~\ref{ReduceOpenToElementary} that $\Xi_{\epsilon,\delta}(U,X)$ holds
if and only if
$$(\forall \bar{\epsilon} < \epsilon)(\exists Y)(Y\mbox{ is elementary and }
\Xi_{\bar{\epsilon},\delta}(Y,X)$$ and from Lemma~\ref{Xi.2} it follows that $\Xi_{\bar{\epsilon},\delta}(Y,X)$
                                                                              is a Boreel statement
because $Y$, being elementary, is closed.
\stopproof
\begin{defin}
If $X$ and $Y$ are\label{WeakXi} subsets of $[0,1]$ then define the relation
$\Xi^*_{\delta,\epsilon}(X,Y)$ to hold if and only if
 for every elementary set $Z$, if
$\lambda(Z) < \epsilon$  then $H^r(X\cap Y \setminus Z) > H^r(Y) - \delta$.
If $X$ and $Y$ are subsets of $[0,1]^{n+1}$ then define the relation
$\Xi^*_{\delta,\epsilon}(X,Y)$ to hold if and only if
$$\Xi^*_{\delta,\epsilon}(\{x\in \pi_1(X) :
\Xi^*_{\delta,\epsilon}(X_x,Y_x)\}, \pi_1(Y)) .$$
\end{defin}

\begin{lemma}
If $X$ and \label{weakXilem} $Y$ are elementary subsets of $[0,1]^n$  then
$\Xi_{\delta,\epsilon}(X,Y)$
holds if and only if $\Xi^*_{{\delta},{\epsilon}}(X,Y)$  holds. \end{lemma}
\proof      One direction is clear. For the other, proceed by induction on $n$. If $n = 1$ then
suppose that $\Xi^*_{{\delta},{\epsilon}}(X,Y)$  holds
and that  $\lambda(Z) < \epsilon$ is such  that $H^r(X\cap Y\setminus Z) \leq H^r(Y) - \delta$.
It will be shown that there is an elementary $Z'$  such that
$\lambda(Z')  < \epsilon$
 and  $H^r(X\setminus Z') \leq H^r(X\setminus Z)$. This clearly suffices.

The existence of $Z'$ will be established by induction on the number of connected components of
$X\cap Y$.
If $X\cap Y = [a,b]$ is an interval then $H^r(X\cap Y\setminus Z)  \geq  (b-a - \lambda(Z) )^r$
because, if $
\{I_i\}_{i\in\omega}$ is a cover of $X\cap Y \setminus Z $ then $\sum_{i=0}^\infty \lambda(I_i) \geq
b-
a - \lambda(Z)$ and hence, since $r < 1$,  $\sum_{i=0}^\infty \lambda(I_i)^r \geq (\sum_{i=0}^\infty
\lambda(I_i) )^r
\geq   (b-a - \lambda(Z) )^r$. Hence $Z'$ can be chosen to be a sufficiently small rational interval
containing $[a, a+\lambda(Z)]$ because then
$H^r(X\cap Y \setminus Z') \leq (b-a - \lambda(Z) )^r  \leq H^r(X\setminus Z)$.
Now suppose that  $X\cap Y  = [a_0,b_0] \cup [a_1,b_1] \cup \ldots \cup [a_k,b_k]$ where $a_i < b_i
< a_{i+1}< b_{i+1}$ for each $ i \in k$. It is possible to
choose open sets $\{U^j_i\}_{i\in\omega}$ such that $X\cap Y \setminus Z\subseteq
\cup_{i\in\omega}U^j_i$
and $\sum_{i=0}^\infty \lambda(U^j_i)^r <   H^r(X \cap pY \setminus Z)  + \frac{1}{j + 1}$ for each
$ j \in \omega$. If  $j \in \omega$ is such that $(b_0,a_{1}) \not\subseteq
        \cup_{i\in\omega}U^j_i$ then  it may as well be assumed that
$U^j_i\cap U^j_{i'} = \emptyset $ if $U^j_i\cap [a_0, b_0] \neq
\emptyset$ and $U^j_{i'}\cap [a_m, b_m] \neq
\emptyset$ for some $m > 0$.
Hence, if there are infinitely many $j \in \omega$ such that $(b_0,a_{1}) \not\subseteq
        \cup_{i\in\omega}U^j_i$ then it
that follows that   $H^r(X\setminus Z)  = $
        $$H^r( [a_0,b_0]\setminus Z) +
         H^r( [a_{1},b_{1}] \cup [a_{2},b_{2}] \cup \ldots \cup [a_k,b_k]\setminus Z)$$ and
the induction hypothesis can be used to find elementary $Z_0$ and $Z_1$
such that
\begin{itemize}
\item  $H^r( [a_0,b_0] \setminus Z_0) \leq
H^r( [a_0,b_0] \setminus Z)$ \item
$H^r( ([a_{1},b_{1}] \cup [a_{2},b_{2}] \cup \ldots \cup
[a_k,b_k])\setminus Z_1) \leq
 H^r(( [a_{1},b_{1}] \cup [a_{2},b_{2}] \cup \ldots \cup
[a_k,b_k])\setminus Z) $
\item
 $\lambda(Z_0)
< \epsilon_0$ \item $\lambda(Z_1) < \epsilon_1$
\end{itemize}
  where $\epsilon_0$ and $\epsilon_1$
are chosen to be positive so that $\epsilon_0+ \epsilon_1 \leq \epsilon$ and
$\lambda([a_0,b_0] \cap Z) < \epsilon_0$
and    $\lambda(([a_{1},b_{1}] \cup [a_{2},b_{2}] \cup \ldots \cup [a_k,b_k])\cap  Z) <\epsilon_1$.

 Hence it may be assumed that for all but finitely many $j\in \omega$
there is some $b(j)\in \omega$ such that
$(b_0,a_{1}) \subseteq
        U^j_{b(j)} = (x_j, y_j)$. By restricting attention to an infinite
subsequence it may also be assumed that there is some interval $[x,y]$
such that $\lim_{j\to\infty}x_j = x$ and  $\lim_{j\to\infty}y_j = y$.
It follows that $$H^r(X\cap Y \setminus Z) = H^r([a_0,x] \setminus Z) +
(y-x)^r + H^r(([y,b_{1}] \cup [a_{2},b_{2}] \cup \ldots \cup
[a_k,b_k])\setminus  Z)$$
and so the induction hypothesis can be used to find
 elementary $Z_0$, $Z_1$ and $Z_2$
such that
\begin{itemize}
\item  $H^r( [a_0,x] \setminus Z_0) \leq
H^r( [a_0,x] \setminus Z)$ \item
$H^r( ([y,b_{1}] \cup [a_{2},b_{2}] \cup \ldots \cup
[a_k,b_k])\setminus Z_1) \leq
 H^r( [y,b_{1}] \cup [a_{2},b_{2}] \cup \ldots \cup [a_k,b_k]\setminus Z) $ \item
 $\lambda(Z_0)< \epsilon_0$
\item $\lambda(Z_1)< \epsilon_1$
\item $Z_2 = J_x\cup J_y$ and $ x \in J_x$ and $ y \in J_y$
\end{itemize}
  where $\epsilon_0$, $\epsilon_1$ and $\epsilon_2$
are chosen to be positive so that $\epsilon_0+ \epsilon_1 + \epsilon_2 \leq \epsilon$
$\lambda([a_0,x] \cap Z) < \epsilon_0$,
  $\lambda(([y,b_{1}] \cup [a_{2},b_{2}] \cup \ldots \cup [a_k,b_k])\cap  Z)
<\epsilon_1$ and
$\lambda([x,y]\cap Z) < \epsilon_2$.
Let $Z = Z_0 \cup Z_1 \cup Z_2$.

Now suppose that the result has been established for $n$ and that $X$
and $Y$ are elementary subsets of $[0,1]^{n+1}$ and that
$\Xi_{\delta,\epsilon}(X,Y)$.
In other words,
$$\Xi^*_{\delta,\epsilon}(\{x\in \pi(X) :
\Xi^*_{\delta,\epsilon}(X_x,Y_x)\}, \pi(Y))$$
holds and, using the induction hypothesis this yields that
$$\Xi^*_{\delta,\epsilon}(\{x\in \pi(X) :
\Xi_{\delta,\epsilon}(X_x,Y_x)\}, \pi(Y))$$
holds as well. To finish use the case $n = 1$.
\stopproof

\begin{lemma} If $K\subseteq [0,1]$ is closed \label{pp30}
then for all $\mu > 0$ there is a closed subset $K' \subseteq K$ such that $\lambda(K\setminus K') <
\mu$ and for each $i\in\omega$ there is $\gamma > 0$ such that $\Xi_{\frac{1}{i+1},\gamma}(K',K')$.
\end{lemma}
\proof First, the following weaker statement will be established:
If $K\subseteq [0,1]$ is closed
then for all $\mu > 0$, $i\in\omega$ there is $\gamma > 0$
 and a closed subset $K' \subseteq K$ such that $\lambda(K\setminus K')
< \mu$ and  $\Xi_{\frac{1}{i+1},\gamma}(K',K')$. To see this, suppose not and that $K$, $\mu$ and $i
\in \omega$ provided a counterexample.  Choose inductively  open sets $A_m$ such that
$\lambda(A_m ) < \mu/2^{i+1} $ and
$$H^r(K \setminus (\cup_{j\leq m} A_j)) \leq H^r(K \setminus (\cup_{j\in  m} A_j))   - 1/i$$
for each $ m\leq i+2$. If it is not possible to do this for some $m$ then it follows that
$\Xi_{\frac{1}{i+1},\mu/2^{i+1}}(K \setminus (\cup_{j\in m} A_j), K
\setminus (\cup_{j\in m} A_j))$ holds and $\lambda(\cup_{j\in m}A_j) < \mu$. On
the other hand, if the induction can be completed then the following inequalities hold:
\begin{itemize}
\item[1] $H^r(K) - \frac{1}{i+1} \geq H^r(K\setminus A_0)$
\item[2]    $H^r(K\setminus A_0) - \frac{1}{i+1} \geq H^r(K\setminus (A_0\cup A_1))$
\item[$\vdots$]
\item[m] $H^r(K \setminus (\cup_{j\in m}A_j) - \frac{1}{i+1} \geq H^r(K \setminus (\cup_{j\leq
m}A_j))$ \end{itemize}
and therefore $H^r(K) - \frac{1}{i+1} \geq  H^r(K \setminus (\cup_{j\leq m}A_j))  \geq 0$
contradicting that $K\subseteq [0,1]$ implies that $H^r(K) \leq 1$.

Now choose inductively open sets $U_i$ and numbers  $\bar{\gamma}(i) > 0$ such that
\begin{itemize}
\item $\bar{\gamma}(0) = \mu$
\item $\lambda(U_i) <  \frac{\bar{\gamma}(j)}{2^{i - j + 2}}$ for each $ j \leq i \in \omega$
\item $\Xi_{\frac{1}{i+1}, \bar{\gamma}(i)}(K \setminus (\cup_{j\leq i}U_j) ,K \setminus (\cup_{j\leq
i}U_j))$  holds.
\end{itemize}
Now let $K' = K \setminus (\cup_{j\in \omega}U_j) $ and note that $\lambda(K\setminus K') < \mu$.
Moreover  $\Xi_{\frac{1}{i+1}, \bar{\gamma}(i)/2}(K' , K')$  holds for each $ i\in \omega$ because
 $\lambda(\cup_{m= i +1}^\omega U_m) < \bar{\gamma}(i)/2$ and hence,
 $\lambda((K \setminus    (\cup_{j\leq i}U_j) )\setminus K'  ) < \bar{\gamma}(i) /2$. Since
 $\Xi_{\frac{1}{i+1}, \bar{\gamma}(i)}(K \setminus (\cup_{j\leq i}U_j) ,K \setminus (\cup_{j\leq
i}U_j))$  holds by construction, the result
                                           follows by setting
$\gamma(i) = \bar{\gamma(i)}/2$.
\stopproof

\section{The Definition of the Ideals Associated with a Capacity}

This section contains the definition of the ideals which will be used
to construct the partial orders satisfying KP$(r)$. Most of the technical
concepts have already been introduced but a few more are needed.

\begin{defin}
A sequence $\{X_i : i\in \omega\}$ of subsets of $[0,1]^d$ will be
said to be a normal family if   \label{normfam}
\begin{itemize}
\item each $X_i$ is elementary
\item $X_{i+1} \subseteq X_i$
\item $\displaystyle \lambda(\pi_n(X_i)\setminus \pi_n(X_{i+1})) <
\frac{\lambda(\pi_n(X_i))}{2^{i+2}}$ for $n\leq d$
\item for each $ i \in \omega$ there is $\beta(i) > 0$ such that $\Xi_{1/i,\beta(i)}(X_j,X_j)$
holds for all $j \geq i$. \end{itemize}
The family
$\{X_i : i\in \omega\}$ will be said to be of dimension $d$. The
function $\beta$ will be called a witness to the normality of the
family $\{X_i : i\in \omega\}$.
\end{defin}
Observe that the intersection of any normal family has positive measure. In fact, if $\{X_i :
i\in \omega\}$ is a normal family and $X = \cap_{i\in \omega}X_i$ then
$$\lambda(\pi_n(X)) > \frac{(2^{i+1} - 1)\lambda(\pi_n(X_i))}{2^{i+1}}$$ for any $i
\in \omega$ and $ n \leq d$.

\begin{defin} Let $W_n$ \label{W} be the family of all sets $a = \bigcup^k_{i=0}I_i\subseteq
[0,1]$ where each $I_i$ is a rational interval and
$\sum_{i=0}^k\lambda(I_i) < 2^{-n}$.  \end{defin}

\begin{defin} Let $n \in \omega$ and $\delta > 0$.  Suppose that $d\in\omega$, ${\cal C} =
\{C_i\}_{i\in\omega}$ is a \label{10} normal family of dimension  $d$ and  $f:[0,1]^d
\rightarrow [0,1]$ is a continuous function.
Define
$X(f,{\cal C},\delta))$ to be the set of all $ a\in W_n$ such that for every $\epsilon > 0$ there
are infinitely many $i\in \omega$ such  that
$\Xi_{\delta,\epsilon}(f^{-1}a,C_i)$ does not hold. The set
$X(f,{\cal C},\delta)$ will be said to be of dimension $d$.
Define ${\cal I}^r_n$ to be the set of all sets $Y\subseteq W_n$ such that
there are $\delta > 0$,   $m \geq 1$,
a normal family ${\cal C}$ of dimension  $ m$   and a continuous function
$f:[0,1]^d \to [0,1]$ such that   $Y \subseteq X(f,{\cal C},\delta)$.  \end{defin}

The ideals of Definition~\ref{10} are defined on the countable set  $W_n$ rather than  on $\omega$.
Theorem~\ref{manu.main} still applies of course.

\begin{lemma}
Let $0 < a < \frac{1}{2^{1/r} + 2}$.
If\label{sum=sum}  $A \subseteq [0,a]$ and $B\subseteq [1- a,1]$ then $H^r(A\cup B) =
H^r(A) + H^r(B)$.
\end{lemma}
\proof Noting that the hypothesis on $a$ implies that $0 < (1 - 2a)^r
- 2a^r$ it is possible to choose $\epsilon > 0$ such that
$\epsilon < (1 - 2a)^r - 2a^r$. Since $H^r(A\cup B) \leq 2a^r$ it
follows that if $A\cup B \subseteq \cup_{i\in\omega}I_i$ and
$\sum_{i\in\omega}\lambda(I_i)^r < H^r(A\cup B) + \epsilon$ then none
of the intervals $I_i$ contains $[a, 1 - a]$. Since $a < 1/2$ it may
as well be assumed that none of the intervals contains $1/2$, or, in
other words, that $\{i \in \omega : I_i \cap A \neq \emptyset\}$ is
disjoint from $\{i \in \omega : I_i \cap B \neq \emptyset\}$. Hence
$H^r(A\cup B) \geq H^r(A) + H^r(B)$. The
result follows since $H^r(A\cup B) \leq H^r(A) + H^r(B)$ is true in general.
\stopproof

\begin{lemma}
Each set ${\cal I}^r_n$ is closed under \label{FiniteUnions} finite unions.
\end{lemma}
\proof
Let $X(f,\{ C_i\}_{i\in\omega},\mu)$ and $X(g,\{
D_i\}_{i\in\omega},\rho)$ be any two generators for ${\cal I}^r_n$
of dimension $d_1$ and $d_2$ respectively.
Let $d $ be greater than $d_1$ and $d_2$ and define
$\bar{C}_i = C_i \times [0,1]^{d - d_1}$ and $\bar{D}_i = D_i \times [0,1]^{d -
d_2}$.
Define $\bar{f}(x_1, x_2, \dots x_d) = f(x_1, x_2, \dots , x_{d_1})$
and $\bar{g}(x_1, x_2, \dots x_d) = g(x_1, x_2, \dots , x_{d_2})$.

Next, let  $0 < a < \frac{1}{2^{1/r} + 2}$.
Define $\psi_1 : [0,1] \to [0,a]$ by $\psi_1(x) = ax$ and define
$\psi_2 : [0,1] \to [1- a, 1]$ by
 $\psi_2(x) = 1 - ax$.
Let $\varphi_i:[0,1]^d \to [0,1]^d$ be defined by $\varphi_i(x_1,x_2,\ldots,
x_d) = (\psi_i(x_1),x_2,\ldots,
x_d)$ for $i \in \{1,2\}$. Let $B_i = \varphi_1(\bar{C}_i) \cup \varphi_2(\bar{D}_i)$.
To see that $\{B_i\}_{i\in\omega}$ is a normal family  it must only be observed that
if $\beta_1:\omega \to (0,1)$ witnesses that $\{ C_i\}_{i\in\omega}$ is normal and
 $\beta_2:\omega \to (0,1)$ witnesses that $\{ D_i\}_{i\in\omega}$ is normal then the function
 $\beta : \omega \to (0,1)$
 defined  by $\beta(i) = \min\{a\beta_1(i) ,a\beta_2(i)\}$ witnesses the normality of
$\{B_i\}_{i\in\omega}$. This uses Lemma~\ref{sum=sum}. Finally, let $h$
be any continuous extension of $(\bar{f}\circ \varphi_1^{-1}) \cup (\bar{g}\circ \varphi_2^{-1})$
and let $\delta = \min\{a^r\mu, a^r\rho\}$. Clearly
$X(h,\{B_i\}_{i\in\omega}, \delta) \in {\cal I}^r_n$.

It will be shown that  $X(f,\{ C_i\}_{i\in\omega},\mu)$ is a subset of
$X(h,\{B_i\}_{i\in\omega}, \delta)$, the proof for $X(g,\{
D_i\}_{i\in\omega},\rho)$ being similar. Let $b \in X(f,\{
C_i\}_{i\in\omega},\mu)$. This means that for every $\epsilon > 0$
there are infinitely many $i\in \omega$ such that
$\Xi_{\mu,\epsilon}((f^{-1}b), C_i)$ fails to hold. Let $\epsilon$ and
$i$ be fixed such that $\Xi_{\mu,\epsilon}((f^{-1}b), C_i)$ fails. Unraveling the
definition of $\Xi_{\mu,\epsilon}$ reveals that
$$H^r(\{x \in \pi_1(C_i) : \Xi_{\mu, \epsilon}((f^{-1}b)_x,
(C_i)_x)\}\setminus Z)
\leq H^r(\pi_1(C_i)) - \mu$$
 for some  set $Z$
such that $\lambda(Z) < \epsilon$. From the definition of $\varphi_1$ and $\bar{f}$ it
follows that
$$H^r(\{x \in \pi_1(\bar{C}_i) : \Xi_{\mu, \epsilon}((h^{-1}b)_{\psi_1(x)},
(B_i)_{\psi_1(x)})\}\setminus Z)
\leq H^r(\pi_1(\bar{C}_i)) - \mu$$
 and so,
observing  that
$H^r(\psi_1(A))= a^rH^r(A)$ for any $A \subseteq
[0,1]$,  it follows that
$$H^r(\{x \in \pi_1(\varphi_1(\bar{C}_i)) : \Xi_{\mu, \epsilon}((h^{-1}b)_{x},
(B_i)_{x})\}\setminus Z') = $$
$$H^r(\{\psi_1(x) : x \in \pi_1(\bar{C}_i) \AND \Xi_{\mu, \epsilon}((h^{-1}b)_{\psi_1(x)},
(B_i)_{\psi_1(x)})\}\setminus Z')$$
$$\leq a^r(H^r(\pi_1(\bar{C}_i)) - \mu) =
H^r(\pi_1(\varphi_1(\bar{C}_i))) - \delta$$
 where $Z'$ is the image of $Z$ under $\psi_1$.
Notice that $\lambda(Z') = a\lambda(Z) < \lambda(Z) < \epsilon$.

The next thing to notice is that
$$\{x \in \pi_1(\varphi_1(\bar{C}_i)) : \Xi_{\mu, \epsilon}((h^{-1}b)_{x},
(B_i)_{x})\}  \supseteq
\{x \in \pi_1(\varphi_1(\bar{C}_i)) : \Xi_{\delta, \epsilon}((h^{-1}b)_{x},
(B_i)_{x})\}$$
because $\delta <
\mu$.
From Lemma~\ref{sum=sum} it follows that $ H^r(\pi_1(B_i)) = H^r(\pi_1(\varphi_1(\bar{C}_i))
+ H^r(\pi_1(\varphi_2(\bar{D}_i))$. Therefore,
$$H^r(\{x \in \pi_1(B_i) : \Xi_{\delta, \epsilon}((h^{-1}b)_{x},
(B_i)_{x})\}\setminus Z')
 \leq$$
$$H^r(\{x \in \pi_1(\varphi_1(\bar{C}_i)) :
\Xi_{\delta, \epsilon}((h^{-1}b)_{x},
(\bar{C}_i)_{x})\}) + H^r(\pi_1(\varphi_2(\bar{D}_i))
\leq $$
$$H^r(\{x \in \pi_1(\varphi_1(\bar{C}_i)) :
\Xi_{\mu, \epsilon}((h^{-1}b)_{x},
(\bar{C}_i)_{x})\}) + H^r(\pi_1(\varphi_2(\bar{D}_i))
\leq $$
$$H^r(\pi_1(\varphi_1(\bar{C}_i))) - \delta +
H^r(\pi_1(\varphi_2(\bar{D}_i)) = H^r(\pi_1(B_i)) - \delta$$
or, in other word
$\Xi_{\delta,\epsilon}(h^{-1}b,B_i)$ fails provided that
$\Xi_{\mu,\epsilon}(f^{-1}b,C_i)$ fails. Since for every $\epsilon >
0$  there are infinitely many $i\in \omega$ such that
$\Xi_{\mu,\epsilon}(f^{-1}b,C_i)$ fails it follows that
 $b \in X(h,\{B_i\}_{i\in\omega},\delta)$.
\stopproof

\begin{lemma} If the parameters $f$, ${\cal C}$ and $\delta$ are given then the statement
\label{Pi11.p} ``$a \in X(f,{\cal C},\delta)$'' is
arithmetic. \end{lemma}
\proof  Let ${\cal C} = \{C_i\}_{i\in\omega}$ be a normal family of dimension  $n$.     From
Definition~\ref{10}
it follows that $a \in X(f,{\cal C},\delta)$ if and only if for every $\epsilon > 0$ there
are infinitely many $i\in \omega$ such  that $\Xi_{\delta,\epsilon}(f^{-1}a,C_i)$ does not hold.
Since $a$ is open and $f$ is continuous it follows from Lemma~\ref{ReduceOpenToElementary}
that $a
\in X(f,{\cal C},\delta)$ if and only if  $$(\forall \epsilon > 0)(\forall m\in \omega)(\exists i >
m)(\exists  \bar{\epsilon} < \epsilon)(\forall Y)( Y\mbox{ is
elementary and }$$
$$ Y\subseteq f^{-1}a
\Rightarrow  \neg \Xi_{\delta,\bar{\epsilon}}(Y,C_i))$$
and $\Xi_{\delta,\bar{\epsilon}}(Y,C_i)$ is equivalent to $\Xi^*_{\delta,\bar{\epsilon}}(Y,C_i)$
when $Y$ and $C_i$ are elementary
by
Lemma~\ref{weakXilem}.

Hence it suffices to show that the statement $\Xi^*_{\delta,\bar{\epsilon}}(Y,C_i)$ is arithmetic.
Proceed by induction on $n$. Notice that the statements $\lambda(Z) < \alpha$  and $H^r(Z) >
\alpha$ are arithmetic for elementary sets $Z$. The case $n=1$ follows immediately
and the induction is carried through because of the elementarity of $Y$ and $C_i$.
\stopproof

\begin{lemma}If the parameters $\beta$ and  ${\cal C}$ are given then the statement
\label{Pi11.b} ``$\beta$ witnesses the normality of ${\cal C}$'' is
arithmetic. \end{lemma}
\proof
This follows from Lemma~\ref{weakXilem} and the definition of a normal
family because it has alraedy been observed in the proof of
Lemma~\ref{Pi11.p} that the statement $\Xi^*_{\delta,\bar{\epsilon}}(Y,C_i)$ is arithmetic.
\stopproof

\begin{corol} The ideals \label{Pi11} $ {\cal I}^r_n$ are all
$\Sigma^1_1$ ideals.
\end{corol}
\proof     From Definition~\ref{10} it follows that $Y \in {\cal I}^r_n$ if and only if
there are $\delta > 0$,   $m \geq 1$,
a normal family ${\cal C}$ of dimension  $ m$   and a continuous function
$f:[0,1]^m \to [0,1]$ such that   $Y \subseteq X(f,{\cal C},\delta)$. Now apply Lemma~\ref{Pi11.p}
noting that the existence of a normal family can be expressed with a $\Sigma^1_1$ statement.
\stopproof

\begin{lemma}
If $A \subseteq B \subseteq [0,1]^d$ and $X \subseteq B$ are such that
\label{inf-to-1}
\begin{itemize}
\item $\lambda(\pi_n(B) \setminus \pi_n(A)) <
(\frac{\epsilon}{d+1})^n$ for each $ n \leq d$
\item $\Xi_{\delta,\epsilon}(X,B)$
\end{itemize}
then $\Xi_{\delta,\frac{\epsilon}{d+1}}(X, A)$.
\end{lemma}
\proof
Proceed by induction in $d$. If $d = 1$ and $\lambda(Z) < \epsilon/2$
then $\lambda((B\setminus A) \cup Z) < \epsilon$.
 Hence $H^r(X \cap A\setminus((B\setminus A) \cup Z)) >
H^r(B) - \delta$. Since $X \cap A\setminus(( B\setminus A)\cup Z)  =
X\cap A \setminus Z$ this suffices.

Suppose the lemma is true for $d$  and that
$A\subseteq B \subseteq [0,1]^{d+1}$. Let
$$S_n = \{x\in [0,1] : \lambda((\pi_n(B_x) \setminus \pi_n(A_x)) \geq
(\frac{\epsilon}{d+2})^n\}$$
for each $ n\leq d$. Since $\lambda((\pi_{n+1}(B) \setminus \pi_{n+1}(A))  <
(\frac{\epsilon}{d+2})^{n+1}$ it follows that $$\lambda(S_n) <
\epsilon/(d+2)$$ for each $ n\leq d$.
If $Z$ is such that $\lambda(Z) < \epsilon/(d+2)$ define
$Y(Z) = Z \cup (\cup_{n=1}^d S_n ) \cup (\pi_1(B) \setminus \pi_1(A))$
and note that
that $\lambda(Y(Z)) < \epsilon$.
Hence
$$H^r(\{x\in \pi_1(B) : \Xi_{\delta,\epsilon}(X_x,B_x)\}\setminus (Y(Z))
> H^r(\pi_1(B)) - \delta$$
and,  moreover, if
$\Xi_{\delta,\epsilon}(X_x,B_x)$ holds and $ x\notin  Y(Z)$ then
$A_x$, $B_x$ and $X_x$ satisfy the hypothesis of the lemma for  $d$
and, furthermore, $x \in \pi_1(A)$.
Therefore
$$H^r(\{x \in \pi_1(A): \Xi_{\delta,\frac{\epsilon}{d+1}}(X_x,A_x)\}\setminus Z
) >
H^r(\pi_1(B)) - \delta >H^r(\pi_1(A)) - \delta $$ and this implies
that
$$H^r(\{x \in \pi_1( A): \Xi_{\delta,\frac{\epsilon}{d+2}}(X_x, A_x)\}\setminus Z
) >H^r(\pi_1(A)) - \delta $$ Since $Z$ was arbitrary, this means that
$\Xi_{\delta,\frac{\epsilon}{d+2}}(X, A)$ holds.
\stopproof

\begin{corol} If $\{C_i\}_{i\in\omega}$ is a normal family of
dimension\label{equivs3} $d$ then  the following are equivalent:
\begin{enumerate}
\item There
is $\epsilon > 0$ such that $\Xi_{\delta,\epsilon} (X,C_i)$ holds for all but finitely many $i \in
\omega$.
\item There
is $\epsilon > 0$ such that $\Xi_{\delta,\epsilon} (X,C_i)$ holds for infinitely many $i \in
\omega$.
\item There
is $\epsilon > 0$ such that $\Xi_{\delta,\epsilon} (X,C_i)$ holds for some $i \in \omega$ such that
$$\lambda(\pi_n(C_i) \setminus \pi_n(\cap _{j\in\omega}C_j)) <
(\frac{\epsilon}{d+1})^n$$ for each $ n \leq d$.
\end{enumerate}      \end{corol}
\proof
To get that  (3) implies (1) use Lemma~\ref{inf-to-1} noting that if $ j > i$ then
$\lambda(\pi_n(C_i) \setminus \pi_n(C_j)) < (\frac{\epsilon}{d+1})^n$ for each $ n \leq d$ and so,
$\Xi_{\delta, \frac{\epsilon}{d+1}}(X,C_j)$ holds.
\stopproof

\begin{lemma} Each of the  ideals \label{6} ${\cal I}^r_n $ of
Definition~\ref{10} satisfies KP$(r)$.  \end{lemma}
\proof Suppose
not.  Then there is some $n\in\omega$ such that the ideal ${\cal
I}^r_n$ does not satisfy KP$(r)$. This means that there exist
\begin{itemize} \item $\theta > 0$
\item $X \in {\cal I}^+$ \item a function $F$ from $X$ to the Borel  subsets of $[0,1]$
satisfying that $H^r(F(x)) \leq \theta$ for each $x\in X$ \item
$\epsilon > 0$
\end{itemize} such that for every $Y\subseteq [0,1]$ and $Z\subseteq [0,1] $ such that
\begin{itemize} \item $H^r(Y) \leq \theta$ \item $\lambda(Z) < \epsilon$
\end{itemize} it must be the case that $\{a \in X :  y \notin F(a)\} \in {\cal I}$ for some $y
\in [0,1] \setminus (Y\cup Z)$.
  Using Definition \ref{10}, it is possible to rephrase this as
follows: For every $Y\subseteq [0,1]$ and $Z\subseteq [0,1] $ such
that $H^r(Y) \leq \theta$ and $\lambda(Z) < \epsilon$ it must be that
there is some $y \in [0,1] \setminus (Y\cup Z)$ and
there are $\delta > 0$,   $m \geq 1$,
a normal family ${\cal C} $ of  dimension  $ m$   and a continuous function
$f:[0,1]^m \to [0,1]$ such that   $\{a \in X : y \notin F(a)\} \subseteq X(f,{\cal
C},\delta)$

Let ${\cal E}_m$ be the set of all elementary subsets of $[0,1]^m$ considered to have the discrete
topology.  It follows that $\prod_{m\in \omega}{\cal E}_m$ is homeomorphic to the irrationals. Let ${\cal
N}_m$ be the subspace of  $\prod_{m\in \omega}{\cal E}_m$ consisting of all $\xi$ such that
$\{\xi(n)\}_{n\in\omega}$ is a normal family and observe that, because it is a closed subspace of
$\prod_{\omega}{\cal E}_m$, ${\cal N}_m$ is a Polish space.
Let ${\cal C}([0,1]^m)$ be the space
of continuous functions from $[0,1]^m$ to $[0,1]$   with the metric induced by the supremum norm.
Let $${\cal P}_m =  {\cal C}([0,1]^m)   \times {\cal N}_m\times (0,1)\times (0,1)^\omega$$
and let ${\cal P} = \cup_{m\in
\omega} {\cal P}_m$ and note
                            that $\cal P$ is a Polish
space.
Let $\Omega$ be the set of all $(z, g, \xi, \delta, \beta) \in [0,1]\times{\cal P}$
such that  $\{a \in X : z\notin F(a)\}\subseteq X(g,\{\xi(n)\}_{n\in\omega},\delta)$ and the
normality of  $\{\xi(n)\}_{n\in\omega}$ is witnessed by $\beta$. Because $X$ and $F$ can be coded by
reals, the definition
of $\Omega$ together with Lemma~\ref{Pi11.p} and Lemma~\ref{Pi11.b} immediately establish that $\Omega$ is a Borel subset
of the Polish space $[0,1]\times{\cal P}$.

It is therefore possible to appeal to the von Neumann Selection
Theorem to find a measurable $\Phi :[0,1] \rightarrow \cal P$ such that the domain of $\Phi$ is the
same as $\pi_1(\Omega)$ and
$\Phi \subseteq
\Omega$.  If $x$ is in the domain of
$\Phi$ suppose that  $\Phi(x) = (g, \xi,\delta,\beta)$    and define $d(x)$
to be  the dimension of $X(g,\{\xi(n)\}_{n\in\omega}, \delta)$.
 Then  define
  $\Phi_i^n(x) = \pi_n(\xi(i))$ for each $i \in \omega$ and
 define  $\Phi^n_\omega(x) =\pi_n(
\cap_{i\in\omega}\xi(i))$ --- if $n
> d(x)$ then $\pi_n(\xi(i)) = \xi(i)$. Since
$\lim_{i\to\infty}\lambda(\Phi_i^n(x)) = \lambda(\Phi^n_\omega(x))$
for each $x $ in the domain of $\Phi$ and $n \in \omega$,
it is possible to apply Egerov's theorem  countably many times to find
 a compact set $\bar{K}$ --- which is the intersection of a nested sequence of closed sets obtained
from the countably many applications of  Egerov's theorem --- such that
\begin{itemize}
\item $\Phi\restriction \bar{K}$ is continuous
\item $\Phi^n_\alpha\restriction \bar{K}$ is continuous  for each $\alpha \in \omega+1$
\item $\{\lambda(\Phi^n_i(x))\}_{i\in\omega}$ converges uniformly, with respect to the variable $x$,
to $\lambda(\Phi^n_\omega(x))$ on $\bar{K}$ \item $\lambda(\bar{K}) > \lambda(   \pi(\Omega)) -
\epsilon/4$. \end{itemize}
 Observe that if
 $Z$ is such that $\lambda(Z) < \epsilon/2$ then $H^r( \bar{K} \setminus Z) > \theta$ because
otherwise, it is possible to obtain a contradiction by setting $Y = \bar{K} \setminus Z$ in the definition of KP$(r)$. Now use
Lemma~\ref{pp30} to find a closed  $K\subseteq \bar{K}$ such that $\lambda(\bar{K} \setminus K) < \epsilon/4$
and there exists $\gamma:\omega \to (0,1)$ such that $\Xi_{\frac{1}{i+1},\gamma(i)}(K,K)$ holds for
all $i\in \omega$.

Next, the compactness of $K$ implies that there is $m\in\omega$ such that $d(x) \in m$ for
each $ x\in K$. Furthermore there is  $\delta
> 0$ such that for every  $x \in K$, if $\Phi(x) = (g, \xi,\delta',\beta)$
then $\delta' > \delta$. Since $H^r(K\setminus Z) > \theta$ for each
$Z\subseteq [0,1]$ such that $\lambda(Z) < \epsilon/4$ it follows that,
   by shrinking $\delta$ if necessary, it may be
assumed that  $H^r(K) > \theta + \delta$. Yet another application
of compactness yields a function $\beta:\omega \to (0,1)$ such that for each $x \in K$, if $\Phi(x)
= (g,\xi,\delta,\beta_x)$ then $\beta_x(i) \geq \beta(i)$  for each $i\in \omega$.

Let $\tau_n = \int_K\lambda(\Phi^n_\omega(x))dx$ for $n \leq m$.
Since $\{\Phi^m_i(x)\}_{i\in\omega}$ is a normal family
for each $x $ in the domain of $\Phi$ it follows from the remarks following
Definition~\ref{normfam} that $$\lambda(\Phi^n_i(x))
< \frac{2^{i+1}\lambda(\Phi^n_\omega(x))}{2^{i+1}-1}$$ for each $ i\in
\omega$ and $n \leq m$. Therefore,
$$\int_K\lambda(\Phi^n_i(x))dx < \frac{2^{i+1}\tau_n}{2^{i+1}-1}$$ and so it is possible to choose an
open  set $L_i$ such that $K\subseteq L_i$ and $$\lambda(L_i\setminus K) +
\int_K\lambda(\Phi^n_i(x))dx                      <
\frac{2^{i+1}\tau_n}{2^{i+1}-1}$$ for each $ n\leq m$ and
$H^r(L_i) < H^r(K) + \frac{1}{2i}$ and $\lambda(L_i \setminus K) < \frac{\gamma(i)}{2}$
 for each $ i\in
\omega$.  Next, using the continuity of $\Phi$ on $K$, choose a family
$\{N_i\}_{i\in\omega}$ such that \begin{itemize}
\item $N_i = [p^i_0,q^i_0] \cup \ldots \cup [p^i_{k(i)},q^i_{k(i)}]$ is elementary for each $i$
\item $K\cap [p^i_j,q^i_j] \neq \emptyset$ for each $i\in \omega$ and $ j \leq k(i)$
\item $K\subseteq N_i \subseteq L_{i+1}$
\item $N_{i+1} \subseteq N_i$
\item $\Phi^n_j(x) = \Phi^n_j(y)$ if $ j \leq i $ and $x$ and $y$
belong to $K$ and  the same component of
$N_i$ \end{itemize}
Let $C_i = \bigcup_{j = 0}^{k(i+2)} [p^{i+2}_j,q^{i+2}_j]\times \Phi^m_{i+2}(z)\times
[0,1]^{m-d(x)}$
for $i\in\omega$ where $z$ is chosen arbitrarily from
$[p^{i+2}_j,q^{i+2}_j]\cap K$
 Then, let $C = \cap_{i\in\omega}C_i$. Observe that
$\lambda(\pi_n(C)) = \tau_n$ for $n \leq m$.

Hence, in order to
show that ${\cal C} =
\{C_i\}_{i\in\omega}$ is a normal family,
first observe that if $j \geq i \geq 1$ and $\lambda(Z) < \gamma(2i - 1)$ then
$$H^r(N_j\setminus Z) \geq H^r(K\setminus Z) \geq  H^r(K) - \frac{1}{2i} \geq
H^r(L_i) - \frac{1}{2i} - \frac{1}{2i}$$ and the last expression is at least as large as $ H^r(N_j)
- \frac{1}{i}$. Hence $\Xi_{\frac{1}{i}, \gamma(2i-1)}(N_j, N_j)$
holds for all $ j \geq i\geq 1$. Now let
$\beta^*(i) = \min\{\gamma(2(i+2)-1), \beta(i+2)\}$. Then
$\Xi_{\frac{1}{i}, \beta^*(i)}(N_j, N_j)$ holds for all $ j \geq i\geq 1$ and so does
$\Xi_{\frac{1}{i}, \beta^*(i)}(\Phi^m_{j+2}(z) , \Phi^m_{j+2}(z))$
 because
$\beta^*(i) \leq \beta(iu+2) \leq \beta_{z}(i+2)$
for
any $z \in K\cap N_{i+2}$. Therefore $\Xi_{\frac{1}{i}, \beta^*(i)}(C_j, C_j)$ holds for all $ j \geq i
$.
  Hence, in order to show that ${\cal C}   $ is a normal family
it suffices
to show that $ \lambda(\pi_n(C_i) \setminus \pi_n(C_{i+1})) <
\frac{\lambda(\pi_n(C_i))}{2^{i+2}}$ for $n\leq m$. To see
this, notice that $$\lambda(\pi_n(C_i) \setminus \pi_n(C_{i+1}))  =
\int_{N_{i+2}}\lambda(\Phi^n_{i+2}(x))dx -\int_{N_{i+3}
}\lambda(\Phi^n_{i+3}(x))dx$$
$$ \leq
\int_{N_{i+2}}\lambda(\Phi^n_{i+2}(x))dx -\int_{K }\lambda(\Phi^n_{\omega}(x))dx
$$
$$\leq\lambda((L_{i+2}\setminus K)) +
\int_K\lambda(\Phi^n_{i+2}(x))dx    -\tau_n                  \leq \frac{2^{i+3}\tau_n}{2^{i+3} - 1}
-\tau_n = \frac{\tau_n}{2^{i+3} - 1} \leq\frac{\lambda(\pi_n(C_i))}{2^{i+2}}$$ for each $i\in\omega$.

Now let $f':C \to [0,1]$ be defined by $f'(x,y) = g(y)$ if $\Phi(x) = (g,\xi,\mu,\zeta)$ and extend $f'$
to a continuous function $f :[0,1]^m \to [0,1]$ arbitrarily. Since $X\notin
{\cal I}^r_n$ there must be some  $a \in X$ such that
$a \notin X(f,{\cal C},\delta)$.  This means that there is some
$\epsilon' > 0$ such that $\Xi_{\delta,\epsilon'}(f^{-1}a,C_i)$
holds  for all but finitely many
 $i\in \omega$.   In particular,
$$H^r(\{x\in N_i : \Xi_{\delta,\epsilon'}((f^{-1}a)_x,(C_i)_x)\}\setminus Z) >
H^r(N_i) - \delta > H^r(K) - \delta >\theta$$ holds for all but finitely many $i \in \omega$
and any $Z$ such that $\lambda(Z) < \epsilon'$.  It may, without loss of generality, be assumed
that $\epsilon' \leq \epsilon/2$.

Using the uniform convergence of $\{\lambda(\Phi^n_i(x))\}_{i\in\omega}$ it is possible to find
$j\in\omega$ such that $ \lambda(\Phi^n_i,x) \setminus \Phi^n_{\omega}(x)) < (\epsilon'/m+1))^n$ for
all $x \in K$, $n \leq m$ and $i > j$.  Let $ i > j$ be such that $\lambda(N_i\setminus K) <
\epsilon'$.
     Since
$H^r(F(a)) \leq \theta$ and
$$H^r(\{x\in  N_i : \Xi_{\delta,\epsilon'}((f^{-1}a)_x,(C_i)_x)\}\setminus (N_i
\setminus K)) >\theta$$ it is possible to choose $y\in K\setminus F(a)$ such that
            $\Xi_{\delta,\epsilon'}((f^{-1}a)_y,(C_i)_y)$ holds.
Observe that if
$\Phi(y) = (g,\xi,\delta',\beta')$ then $\xi(n) = (C_n)_y$ for $ n >0$, $g =
f_y$ and $\delta < \delta'$.  The choice of  $j$ guarantees that the
hypothesis (3) of
Corollary~\ref{equivs3} is satisfied by
$i$, $\delta$, $\epsilon'$, $(f^{-1}a)_y$ and $\{(C_n)_y\}_{n\in\omega}$.
 It follows that  there is some $\epsilon > 0$ such that
 $\Xi_{\delta,\epsilon}((f^{-1}a)_y, (C_i)_y)$ holds for all but finitely many $i \in \omega$ and
hence so does
$\Xi_{\delta',\epsilon}((f^{-1}a)_y, (C_i)_y)$.
 Therefore $a \notin
X(g,\{\xi(n)\}_{n\in\omega},\delta')$.
This yields a contradiction to the fact that  $y
\notin F(a)$ and $\Phi(y) = (g,\xi,\delta',\beta')$ implies that $a \in
X(g,\{\xi(n)\}_{n\in\omega},\delta')$.
 \stopproof

\section{The Ideal is Proper}
It remains to be shown that the ideals ${\cal I}^r_n$ are proper. This
will require a careful analysis of the capacity $H^r$. This will require some generalizations of
results from \cite{step.33}. The  key fact about Hausdorff capacity that
will be used is that if $B\subseteq E$ is of small Lebesgue
measure but evenly distributed throughout $E$, then $H^r(B)$ will be
close to $H^r(E)$. This is made precise in the next lemma whose statement
requires the following notation.
\begin{defin}
For any measurable set $A\subseteq [0,1]$ define $\Delta^i_m(A)$ to be
the least real number such that $\lambda(A\cap [0,\Delta^i_m(A)]) =
\frac{i\lambda(A)}{m}$.
\end{defin}
Notice that $\Delta^i_m(A) $ is always defined and that if $A = [0,1]$
then $\Delta^i_m(A) $ is nothing more than $ \frac{i}{m}$.

\begin{lemma} Let  \label{jaut}
$\delta > 0 $,  $ \eta > 0$ and suppose that $E\subseteq [0,1]$ is measurable.
If $\Xi_{\delta,\eta}(E,E)$ holds and $\delta < H^r(E)$ then there
exists  $m \in \omega$
such that ff  $D\subseteq E$ is any measurable
set such that for each $i \in m$ $$ \lambda(D\cap [\Delta^i_m(E),
\Delta^{i+1}_m(E)]) \geq
\frac{\eta}{m}$$ then
$\Xi_{\delta, \frac{\eta}{2m}}(D, E)$.
\end{lemma}
\proof Let $m \in \omega$ be so large that the inequality
$$\frac{ m^{1-r}\eta^{1  + r}}{8\cdot 2^r} > 1$$ is satisfied.
To begin, note that Lemma~\ref{ExTech} implies that
there exists $\bar{\epsilon} > 0$ such that
$\Xi_{\delta -
\bar{\epsilon}, \eta/2}(E,E)$ holds. If $\Xi_{\delta, \frac{\eta}{2m}}(D, E)$
fails then there is some $Z$ such that $\lambda(Z) < \frac{\eta}{2m}$
and an open cover $D \setminus Z \subseteq \bigcup_{i = 0}^\infty I_i$ such that
$\sum_{i=0}^\infty\lambda(I_i)^r < H^r(E) - (\delta - \bar{\epsilon})$.  Let $B = \{i
\in \omega : \lambda(I_i) \geq \frac{1}{2m}\}$ and let $C = \{i \in m
: (\forall j \in B)(I_j\cap [\Delta^i_m(E), \Delta^{i+1}_m(E)] \cap E
= \emptyset\}$.  Three separate cases, depending on the size of $B$
and $C$, will be considered.
\medskip

{\noindent{\bf Case 1}} To begin, suppose that $|B| \geq \frac{
m\eta}{8}$.  Then $$\sum_{i=1}^\infty\lambda(I_i)^r \geq \sum_{i\in
B}\lambda(I_i)^r \geq |B|(1/2m)^r \geq
\frac{ m^{1-r}\mu}{8\cdot 2^r} > 1$$ Since
$\sum_{i=0}^\infty\lambda(I_i)^r < H^r(E) - (\delta - \eta/2) < 1 $ this is
impossible.
\medskip

{\noindent{\bf Case 2}} Suppose now that $|B| < \frac{m \eta}{8}$ and
$|C| \leq \frac{m\eta }{4}$.  It then follows that if $$G= \{ i \in m :
[\Delta^i_m(E), \Delta^{i+1}_m(E)] \cap E \not\subseteq
\bigcup_{j=1}^\infty I_j\}$$ then $|G| \leq 2\cdot |B| + |C|$. The reason for this is that if $j\in
B$ then there are at most two integers $i$ such that the intervals
$[\Delta^i_m(E), \Delta^{i+1}_m(E)]$ intersect $I_j$ but are not
contained in $I_j$ --- this accounts for the summand $2\cdot |B|$.
All the other intervals $[\Delta^i_m(E),\Delta^{i+1}_m(E)]$ for $i \in
G$ must be disjoint from $I_j$ for every $j \in B$ --- this accounts
for the other summand $|C|$.

By the assumptions of this case it follows that $2 \cdot |B| + |C| <
m\eta/2 $ and hence $$\lambda(\bigcup_{i\in G} [\Delta^i_m(E),
\Delta^{i+1}_m(E)]\cap E) < \eta/2$$ Since $\Xi_{\delta -
\bar{\epsilon}, \eta/2}(E,E)$ holds it may be concluded that
$H^r(E \setminus \bigcup_{i\in G} [\Delta^i_m(E),
\Delta^{i+1}_m(E)]) > H^r(E) - (\delta - \bar{\epsilon})$.   Since $E\setminus(\bigcup_{i\in G}
[\Delta^i_m(E), \Delta^{i+1}_m(E)]) \subseteq
\bigcup_{i=1}^\infty I_i$ this yields a contradiction.
\medskip

{\noindent{\bf Case 3}} Suppose that $|B| < \frac{m\eta}{8}$ and $|C| >
\frac{ m\eta}{4}$.
  Let $C'$ be a family of non-consecutive members of $C$ of maximal
cardinality --- hence, $|C'|
\geq |C|/2 >\frac{m\eta}{8}$.  Let $$U_j = \{i \in \omega : I_i\cap [\Delta^j_m(E), \Delta^{j+1}_m(E)]
\cap E  \neq \emptyset\}$$ for each $j
\in C'$ and define $U = \cup_{j\in C'}U_j$.  Since, for $j \in C$, the sets
$[\Delta^j_m(E), \Delta^{j+1}_m(E)]\cap E $ are intersected only by
intervals $I_i$ where $i \in
\omega\setminus B$, and such intervals $I_i$ are smaller than any $[
\Delta^{j+1}_m(E), \Delta^{j}_m(E)]$, it follows that $U_j\cap U_k = \emptyset$ if $k$ and $j$
are distinct members of $C'$.  Therefore, using the fact that $0 < r <
1$, $$\sum_{i\in U} \lambda(I_i)^r \geq \sum_{j\in C'}\sum_{i \in
U_j}\lambda(I_i)^r \geq $$ $$ \sum_{j\in C'}\left(\sum_{i \in
U_j}\lambda(I_i)\right )^r \geq \sum_{j\in
C'}\lambda(D\cap[\Delta^j_m(E),
\Delta^{j+1}_m(E)])^r
\geq$$ $$ \sum_{j\in C'}(\frac{\mu}{m} - \lambda(Z))^r \geq \frac{ m\mu}{8}( \frac{\mu}{2m})^r
 > 1 $$ and once again, as in the first case, this is a contradiction
because $D\subseteq [0,1]$.  \stopproof

If $X\subseteq [0,1]$ then $F:X\to [0,1]$ will be said to have small
fibres if and only if $\lambda(F^{-1}\{x\}) = 0$ for each $ x\in
[0,1]$.  The proof of the Theorem~\ref{main} and the lemmas preceding
it will rely on decomposing an arbitrary continuous function into a
piece that has small fibres and a piece which has countable range.

\begin{lemma} Let $\mu \in (0,1)$ and suppose that $\{X_i : i\in \omega\}$
is a sequence of mutually independent $\{0,1\}$-valued random
variables
 with mean $\mu$ for each $ i\in \omega$. Suppose \label{manylp} that
$C\subseteq [0,1]$ is a measurable set and that for each $j \in n$ the
function $F_j: C \to [0,1]$ is measurable with small fibres. For any
$\rho > 0$ there is $M \in \omega$ such that for all $m \geq M$
the probability that $$\lambda\left (\bigcap_{j\in n}\bigcup_{i \in m}
F_j^{-1}[\frac{i}{m},\frac{i + X_i}{m}]\right ) > \frac{\mu^n\lambda(C)}{2}$$
is greater than $1 - \rho$.
\end{lemma}
\proof  To begin, let $m\in\omega$ be fixed.
For any function $\xi \in {}^nm$ define $\theta(\xi) = \lambda(\bigcap_{j\in
n}F_j^{-1}[\frac{\xi(j)}{m}, \frac{\xi(j) + 1}{m}])$ and let $$Y(\xi)
= \lambda\left (\bigcap_{j\in n}F_j^{-1}[\frac{\xi(j)}{m}, \frac{\xi(j) +
X_{\xi(j)}}{m}]\right ) =
\theta(\xi)\prod_{j\in n}X_{\xi(j)}$$ If $\xi
\neq \xi'$ then
$$\lambda\left (\left (\bigcap_{j\in n}F_j^{-1}[\frac{\xi(j)}{m}, \frac{\xi(j) +
1}{m}]\right )\cap\left (\bigcap_{j\in n}F_j^{-1}[\frac{\xi(j)}{m}, \frac{\xi(j) +
1}{m}]\right )\right  ) = 0$$ and so $\sum_{\xi  \in {}^nm}\theta(\xi) = \lambda(C)$.

Letting $E[Z]$ denote the average value of the random variable $Z$ and
$V[Z]$ the variance of $Z$, it is easy to see that $E[Y(\xi)] =
\theta(\xi)\mu^{\sigma(\xi)}$ where $\sigma(\xi)$ represents the
cardinality of the range of $\xi$. Noting that $\sigma(\xi) \leq n$
for all $\sigma$, it follows that $E[\sum_{\xi \in {}^nm } Y(\xi)] \geq
\mu^n\lambda(C)$. Furthermore,
$$V[\sum_{\xi  \in {}^nm} Y(\xi)] = E[(\sum_{\xi  \in {}^nm} Y(\xi) - E[Y(\xi)])^2]
= $$ $$\sum_{\xi  \in {}^nm}\sum_{\xi'  \in {}^nm}E[(Y(\xi) -
E[Y(\xi)])(Y(\xi') - E[Y(\xi')])].$$ If $\xi$ and $\xi'$ have disjoint
ranges then $Y(\xi)$ and $Y(\xi')$ are independent random variables
and so $$E[(Y(\xi) - E[Y(\xi)])(Y(\xi') - E[Y(\xi')])] = E[Y(\xi) -
E[Y(\xi)]]E[Y(\xi') - E[Y(\xi')]] = 0$$ while if the ranges of $\xi$ and
$\xi'$ are not disjoint then
$$E[(Y(\xi) - E[Y(\xi)])(Y(\xi')
-E[Y(\xi')])] =$$
$$ E[Y(\xi)Y(\xi')] - E[E[Y(\xi)]Y(\xi')] -
E[E[Y(\xi')]Y(\xi)] + E[Y(\xi)]E[Y(\xi')]]$$
$$ =E[Y(\xi)Y(\xi')] -
E[Y(\xi)]E[Y(\xi')] \leq E[Y(\xi)Y(\xi')] $$
$$ =   E[\theta(\xi)\prod_{j\in
n}X_{\xi(j)}\theta(\xi')\prod_{j\in n}X_{\xi'(j)}] \leq \theta(\xi)\theta(\xi')$$
since $X_i \in \{0,1\}$ for each $ i$. It
may be concluded that $$V[\sum_{\xi  \in {}^nm} Y(\xi)]
\leq \sum_{j\in n}\sum_{\xi \in {}^nm }\sum_{\stackrel{\xi'  \in {}^nm}{
\xi'(j) = \xi(j)}}\theta(\xi)\theta(\xi') =
\sum_{j\in n}\sum_{\xi  \in {}^nm}\theta(\xi)\sum_{\stackrel{\xi' \in {}^nm }{
\xi'(j) = \xi(j)}}\theta(\xi') .$$

However, if $j$ is fixed then $\sum_{\stackrel{\xi'  \in {}^nm}{
\xi'(j) = \xi(j)}}\theta(\xi') = \lambda(F_j^{-1}[\frac{\xi(j)}{m} ,
\frac{\xi(j) + 1}{m}])$. Therefore, all that needs to be done is to
choose $M$ so large that if $m \geq M$ and $ i\in m$ then
$\lambda(F_j^{-1}[\frac{i}{m} ,
\frac{i + 1}{m}]) < \frac{\rho\mu^{2n}\lambda(C)}{4n}$ for
each $ j \in n$. The reason this suffices is that this implies that
$$V[\sum_{\xi  \in {}^nm} Y(\xi)]
\leq \sum_{j\in n}\sum_{\xi  \in {}^nm}\theta(\xi)\frac{\rho\mu^{2n}\lambda(C)}{4n}
\leq \frac{\rho\mu^{2n}\lambda(C)^2}{4}$$
and so Chebyshev's Inequality can be applied to conclude that the
probability that $$|\sum_{\xi  \in {}^nm} Y(\xi) - E[\sum_{\xi \in
{}^nm }
Y(\xi)]| > \frac{\mu^{n}\lambda(C)}{2}$$ is less than $\rho$.
Since it has already been established that $E[\sum_{\xi \in {}^nm }
Y(\xi)] \geq
\mu^n\lambda(C)$ it follows that the probability that
$\sum_{\xi  \in {}^nm} Y(\xi) \geq
\mu^n\lambda(C)/2$ is at least $1 - \rho$ as required.

To choose $M$ so large that if $m \geq M$ and $ i\in m$ then
$\lambda(F_j^{-1}[\frac{i}{m} ,
\frac{i + 1}{m}]) < \frac{\rho\mu^{2n}\lambda(C)}{4n}$ for
each $ j \in n$, all that is required is compactness and the fact that
each $F_j$ has small fibres. Since $F_j^{-1}\{x\} =
\cap_{k\in\omega}F_j^{-1}[x-1/k, x + 1/k]$ and
$\lambda(F_j^{-1}\{x\}) = 0$ it follows that it is possible to choose a
finite cover of $[0,1]$ by open intervals, $\cal C$, such that if $I
\in \cal C$ then $\lambda(F_j^{-1}I) <
\frac{\rho\mu^{2n}\lambda(C)}{4n}$ for each $j\in n$. Hence $M$
must be  chosen so large that  if $m\geq M$ and $ i \in m$ then there is
$I\in \cal C$ such that $[\frac{i}{m},\frac{i+1}{m}]\subseteq I$.
\stopproof

\begin{lemma}
Suppose that $\delta > 0 $,  $\mu > 0$, $\eta > 0$ and $k \in
\omega$.  There is then a real number
$\epsilon(\delta,\mu, \eta, k) > 0$ such \label{ions}that if
\begin{itemize}
\item $\{C_i\}_{i \in
k}$ is a family of measurable subsets of $[0,1]$
\item $F_i: C_i \to [0,1]$ is a measurable functions with small fibres
for each $i \in
k$
\item $E \subseteq [0,1]$ is a measurable set
\item $\Xi_{\delta, \eta}(E,E)$ holds and $\delta < H^r(E)$
\item  $\rho > 0$
\end{itemize}
then there
is $M\in \omega$ such that for all $ m > M$ and for
any mutally independent, $\{0,1\}$-valued random variables
$\{X_i\}_{i\in m}$ with mean $\mu$, the probability that
$$\Xi_{\delta,\epsilon(\delta,\mu,\eta,k)}\left (\bigcap_{i\in k}(F^{-1}_i \bigcup_{j\in
m}[\frac{j}{m},\frac{j+X_j}{m}])\cup ([0,1]\setminus C_i), E\right )$$ holds
 is greater than $1 -
\rho$.

Moreover, there is  $\theta > 0$ such that
if \begin{itemize}
\item $E'\subseteq [0,1]$ is a measurable set such that
$\lambda(E\Delta E') < \theta$
\item $\{C_i'\}_{i\in k}$ is a family of measurable sets such that
$\lambda(C_i \Delta C_i') < \theta$ for each $ i\in k$ \item
$\{G_i\}_{i\in k}$ is a family of measurable
functions such that
$$\sup\{|F_i(x) - F_i'(x)| : x\in C\cap C'\} < \theta$$ for each $ i\in
k$
\end{itemize}
then the
probability that
$$\Xi_{\delta,\epsilon(\delta,\mu,\eta,k)}\left (\bigcap_{i\in k}(G^{-1}_i \bigcup_{j\in
m}[\frac{j}{m},\frac{j+X_j}{m}])\cup ([0,1]\setminus C'_i), E'\right )$$ holds
 is still greater than $1 -
\rho$.
 \end{lemma}
\proof Let $\alpha = \mu^k/2$ and use
Lemma~\ref{jaut} to find $p$ such that if $D\subseteq E$ is a
measurable set such that for each $i \in p$ $$ \lambda(E\cap
[\Delta^{i}_{p}(E), \Delta^{i+1}_{p}(E)]) \geq
\frac{\alpha}{2p}$$ then
$\Xi_{\delta, \frac{\alpha}{4p}}(D, E)$ holds. Let $\epsilon(\delta,\mu,\eta,k) =  \frac{\alpha}{4p}$.
Let $\{P_i : i \in s\}$ enumerate the sets of positive measure which
belong to  the coarsest partition of $E$ refining each of the partitions
$\{[\Delta^{i}_{p}(E),\Delta^{i+1}_{p}(E)]\cap E : i \in p\}$ and $\{C_i\cap E, E \setminus
C_i\}$ for $i\in k$. Now use
Lemma~\ref{manylp} to find
 $M \in \omega$
such that for all $m \geq M$ the probability that
$$\lambda\left (\bigcap_{j\in k}\left(( F_j\restriction P_n)^{-1}\bigcup_{i \in m}
 [\frac{i}{m},\frac{i +
X_i}{m}] \cup (P_n \setminus C_j)\right )\right) > \alpha\lambda(P_n)$$ is greater than $1 -
\frac{\rho}{s}$
for each $n \in s$.

Now notice that if $ m \geq M$ is fixed then, because each $F_i$ has
small fibres, it is possible
to find $p(i,j)$ and $q(i,j)$ such that  $ \frac{i}{m} < p(i,j) <
                     q(i,j) < \frac{i+1}{m}$ and
$$\lambda(F_j^{-1}[\frac{i}{m}, \frac{i+1}{m}]) -
\lambda(F_j^{-1}[p(i,j),q(i,j)]) < \frac{\alpha\lambda(P_n)}{2(2k+1)m}$$ for
each $j \in k$ and $ i \in m$. Now observe that if $\lambda(C_i\Delta
C_i') <\frac{\alpha\lambda(P_n)}{2(2k+1)}$ for each $i\in k$ and if
$\lambda(E\Delta E') < \frac{\alpha^\lambda(P_n)}{2(2k+1)}$  and if
$Y_i \in \{0,1\}$ are such that
$$\lambda\left (\bigcap_{j\in k}\left(( F_j\restriction P_n)^{-1}\bigcup_{i \in m}
 [\frac{i}{m},\frac{i +
Y_i}{m}] \cup (P_n \setminus C_j)\right )\right) >
\alpha\lambda(P_n)$$ then the Lebesgue measure of
$$\bigcap_{j\in k}\left(( F_j\restriction P_n)^{-1}\bigcup_{i \in m}
 [p(i,j), p(i,j) + Y_i(q(i,j) - p(i,j))]
 \cup (P_n \setminus C_j)\right )$$ is greater than 
$$\frac{\alpha\lambda(P_n)}{2}.$$
 Therefore, if $\theta > 0$ is such that \begin{itemize}
\item $\theta < \frac{\alpha\lambda(P_n)}{2(2k+1)}$ for each $n\in s$
\item $\theta < p(i,j) - \frac{i}{m}$  for all $i$ and $j$
\item $\theta < \frac{i_1}{m} - q(i,j)$  for all $i$ and $j$
\end{itemize}
then if $\{G_i\}_{i\in k}$ is a family of measurable
functions such that
$\sup\{|F_i(x) - F_i'(x)| : x\in C\cap C'\} < \theta$ for each $ i\in
k$ then $F_j^{-1}[\frac{i}{m},\frac{i + 1}{m}] \subseteq
G_j^{-1}[p(i,j), q(i,j)]$ and hence
$$\lambda\left (\bigcap_{j\in k}\left(( G_j\restriction P_n)^{-1}\bigcup_{i \in m}
 [\frac{i}{m},\frac{i +
Y_i}{m}] \cup (P_n \setminus C'_j)\right )\right) >
\frac{\alpha\lambda(P_n)}{2}$$
also holds for each $ n \in s$.

It follows that the Lebesgue measure of the intersection of the interval
$[\Delta^i_p(E'),\Delta^{i+1}_p(E')]$ with $$ \bigcap_{j\in k}\left(( G_j\restriction P_n)^{-1}\bigcup_{i \in m}
 [\frac{i}{m},\frac{i +
Y_i}{m}] \cup ([\Delta^i_p(E'),\Delta^{i+1}_p(E')] \setminus
C'_j)\right )$$ is greater than
$$\frac{\alpha\lambda([\Delta^i_p(E'),\Delta^{i+1}_p(E')]\cap E')}{2}$$ and the
result  now follows from Lemma~\ref{jaut}.
\stopproof

\begin{lemma}            Let $k \in \omega$ and $\{C_i\}_{i\in k}$ be a family of measurable subsets
of $[0,1]$.
Let $F_i: C_i \to [0,1]$ be a measurable function for each $i \in k$.
 \label{finitelymanyfunctions}
Suppose also that $ \delta > 0$ and $\eta > 0$.
 Then, for any $N\in \omega$
 and $\epsilon > 0$, if $\Xi_{\delta, \epsilon}(E,E)$ holds for some
measurable set $E$ such that $H^r(E) > \delta$ then
 $$\Xi_{\delta, \epsilon} (\bigcap_{i\in k}(F_i^{-1}a)\cup
([0,1]\setminus C_i), E)$$ holds for some $a \in W_N$.
\end{lemma} \proof   For each $i \in k$ let $\{y_j^i : j\in d_i
\leq
\omega\}$ enumerate all points $y \in [0,1]$ such that $\lambda(F^{-1}_i\{y\}) > 0$. Let $C'_i = C_i
\setminus F^{-1}_i\{y^i_j : j\in d_i\}$ and let $F'_i = F_i\restriction C'_i$. Since $F'_i$ has
small fibres for each $i\in k$ it follows from Lemma~\ref{ions} that it is possible to choose $m$ so
large that if $\{X_i\}_{i\in m}$ are $\{0,1\}$-valued random variables
with mean  $2^{-N-1}$ then,
letting $\epsilon' = \epsilon(\delta,2^{-N-1},\eta, k)$, the
probability that  $$\Xi_{\delta, \epsilon'}\bigcap_{i\in
k}((F_i')^{-1} \bigcup_{j\in
m}[\frac{j}{m},\frac{j+X_j}{m}])\cup ([0,1]\setminus C_i'), E)$$ holds
is at least $3/4$  for any measurable set $E$ such that $\Xi_{\delta, \epsilon}(E,E)$ holds
 $H^r(E) > \delta$.  Since the mean of each $X_i$ is $2^{-N-1}$ it is
possible to choose $m$ so large that the probability that
 $$\lambda( \bigcup_{j\in m}[\frac{j}{m},\frac{j+X_j}{m}])) < 2^{-N}$$
is also greater than $3/4$.
Hence, given $E$ with the required properties,
there is $a_0 \in W_N$
  such that $$\Xi_{\delta, \epsilon'}\bigcap_{i\in
k}((F_i')^{-1} a_0)\cup ([0,1]\setminus C_i'), E)$$ holds. Now choose
$J\in \omega$ such that $\lambda( \cup_{i\in k}\cup_{j \geq
J}F_i^{-1}\{y_j^i\}) < \epsilon'/2$.
                It is then easy to find $a \in W_N$ be such that $a_0 \cup \{y^i_j : i
\in k, j \in J\}\subseteq a$. Let $\epsilon = \epsilon'/2$ and note
that it follows that $$\Xi_{\delta, \epsilon}\bigcap_{i\in
k}((F_i)^{-1} a)\cup ([0,1]\setminus C_i), E)$$ holds because, if
$$ Y = \bigcap_{i\in
k}((F_i')^{-1} a_0)\cup ([0,1]\setminus C_i')\setminus \bigcap_{i\in
k}((F_i)^{-1} a)\cup ([0,1]\setminus C_i)$$ then
$Y \subseteq \cup_{i\in k}\cup_{j\geq J}F_i^{-1}\{y^i_j\}$ and hence
$\lambda(Y) < \epsilon'/2$.
\stopproof

\begin{lemma}
Suppose \label{average} that                 $k\in \omega$ and
$\{C_i\}_{i\in k}$ are measurable subsets of $ [0,1]^{d+1}$
and  $F_i:C_i \to [0,1]$ are measurable functions
such that $(F_i)_x$ has small fibres for each  $ x\in [0,1]^d$. Let
$N\in\omega$, $\delta > 0$ and $\eta > 0$.
 Then there is $\epsilon > 0$
such that for all closed $E \subseteq
[0,1]^{d+1}$
and  $\rho > 0$ there is some  $a\in W_n$ such that
the Lebesgue measure of
$$\{x\in \pi_d(E) : \Xi_{\delta,{\epsilon}}
(\bigcap_{i\in k}((F_i^{-1}a)\cup ([0,1]^{d+1}\setminus C_i))_x, E_x) \mbox{ or
}\neg \Xi_{\delta, \eta}(E_x,E_x)\mbox{ or
} H^r(E_x) \leq \delta\}$$
is at least $\lambda(\pi_d(E))(1  - \rho)$.
\end{lemma}
\proof   Let $\{X_i\}_{i\in \omega}$ be a sequence of mutually independent random variables with mean
$2^{-N-1}$. Let $\epsilon = \epsilon(\delta, 2^{-N-1}, \eta, k)$ and
suppose that $E\subseteq [0,1]^{d+1}$ is closed. Next, choose compact subsets $W_i \subseteq \pi_d(C_i)$ and $V_i
\subseteq [0,1]^{d}\setminus W_i$ for $i \in k$ as well as
$E'\subseteq \pi_d(E)$ such that
\begin{itemize}
\item $\lambda(\pi_d(E) \setminus E') <
\frac{\rho\lambda(\pi_d(E))}{6(1 - \rho/2)}$
\item $\lambda([0,1]^d \setminus (\cap _{i\in k}(W_i\cup V_i)) ) < \frac{\rho\lambda(\pi_d(E))}{6(1 - \rho/2)}$
\item  $F_i \restriction W_i$ is continuous
\item  the mapping from $\pi_d(E')$  to $[0,1]$ defined by
$x \mapsto \lambda(E_x\cap (p,q))$ is continuous for
 each pair of rationals $p$ and $q$ such that $0 \leq p < q \leq 1$
\item  the mapping from $\pi_d(W_i)$  to $[0,1]$ defined by
$x \mapsto \lambda((C_i)_x\cap (p,q))$ is continuous for each $ i \in
k$ and each pair of rationals $p$ and $q$ such that $0 \leq p < q \leq 1$ .
\end{itemize}
An easy application of the Lebesgue Density Theorem shows that a
  consequence of the last clause is that if $x\in \pi_d(W_i)$ then
$\lim_{y\to x}\lambda(W_y\Delta W_x) = 0$.  The penultimate clause
implies a similar assertion for $E'$. It is possible to find compact  $E_1$ and $E_2$,
subsets of $E'$ such that
\begin{itemize}
\item if $x \in E_1$ then $\Xi_{\delta, \eta}(E_x,E_x)$ holds and
$H^r(E_x) > \delta$
\item if $x \in E_2$ then $\Xi_{\delta, \eta}(E_x,E_x)$ fails or $H^r(E_x) \leq \delta$
\item $\lambda(E'\setminus (E_1\cup E_2)) <\frac{\rho\lambda(\pi_d(E))}{6(1 - \rho/2)} $
\end{itemize}
because Lemma~\ref{Xi.2} implies that $\{ x \in E' : \Xi_{\delta,
\eta}(E_x,E_x)\}$ is measurable.
 Let $Z = E_1\cap (\cap
_{i\in k}(W_i\cup V_i))$ and for
   $x \in Z$ let $K(x) = \{i\in k : x\in \pi_d(W_i)$ and
notice that $K(x)$ is constant on a neighbourhood of $x$ because the
sets $V_i$ and $W_i$ are all compact.
  If $x \in Z$ and  $i \in K(x)$ then  $F_i\restriction (W_i)_x$
has small fibres, $H^r(E_x) > \delta$ and $\Xi_{\delta, \eta}(E_x,E_x)$ holds,
 so it follows from Lemma~\ref{ions}  that
there is  $\theta_x > 0$ and $M_x\in \omega$ such that
if $||y - x|| < \theta_x$ and $ M \geq M_x$ then
the
probability that
$$\Xi_{\delta,\epsilon}\left ((\bigcap_{i\in K(y)}(F^{-1}_i \bigcup_{j\in
m}[\frac{j}{m},\frac{j+X_j}{m}])\cup ([0,1]\setminus C_i))_y, E_y\right )$$ holds
 is greater than $1 -
\rho^2/2$.

 Since $Z$ is compact,  it is possible to find a single $M$ such that
for all $ m > M$ and  for any $x\in Z$  the probability that
$$\Xi_{\delta,\epsilon}\left ((\bigcap_{i\in K(x)}(F^{-1}_i \bigcup_{j\in
m}[\frac{j}{m},\frac{j+X_j}{m}])\cup ([0,1]\setminus C_i))_y, E_y\right )$$ holds
 is greater than $1 -
\rho^2/2$

Now let $m > M$ be so great that the probability that
$\lambda(\cup_{j\in m}[\frac{j}{m},\frac{j + X_j}{m}]) < 2^{-N}$ is greater than $
2\rho$. Define
$$\Gamma(X_0,X_1,\ldots,X_m)$$ to be the Lebesgue measure of
the set of all $x\in Z$ such that
$$\Xi_{\delta,\epsilon}\left ((\bigcap_{i\in K(x)}(F^{-1}_i \bigcup_{j\in
m}[\frac{j}{m},\frac{j+X_j}{m}])\cup ([0,1]\setminus C_i))_x, E_x\right )$$ holds.
Note that  Corollary~\ref{Xi.2.c} implies that this set is measurable.
 The first step is to estimate $$\alpha_m = \sum_{X_0 = 0}^1 \sum_{X_1 = 0}^1\ldots \sum_{X_{m}=
0}^1 \Gamma(X_0,X_1,\ldots,X_m) \prod_{i=0}^m \mu^{X_i}(1-\mu)^{1 - X_i} $$
the average value of
$\Gamma(X_0,X_1,\ldots,X_m)$. To this end,
let $$\Lambda_x(X_0,X_1,\ldots,X_m) \in \{0,1\}$$ be defined to be 1 if and only if
$$\Xi_{\delta,\epsilon}\left ((\bigcap_{i\in K(x)}(F^{-1}_i \bigcup_{j\in
m}[\frac{j}{m},\frac{j+X_j}{m}])\cup ([0,1]\setminus C_i))_x, E_x\right )$$ holds.
and observe that $\alpha_m  $ is equal to
            $$  \sum_{X_0 = 0}^1 \sum_{X_1 = 0}^1\ldots \sum_{X_{m}=
0}^1\left( \int_{x\in Z}
\Lambda_x(X_0,X_1,\ldots,X_m) dx \right)  \prod_{i=0}^m \mu^{X_i}(1-\mu)^{1 - X_i} = $$
$$  \int_{x\in Z} \left( \sum_{X_0 = 0}^1 \sum_{X_1 = 0}^1\ldots \sum_{X_{m}=
0}^1 \Lambda_x(X_0,X_1,\ldots,X_m)
            \prod_{i=0}^m \mu^{X_i}(1-\mu)^{1 - X_i} \right)dx  $$
 However, notice that
   $$ \sum_{X_0 = 0}^1 \sum_{X_1 = 0}^1\ldots \sum_{X_{m}=
0}^1 \Lambda_x(X_0,X_1,\ldots,X_m)
            \prod_{i=0}^m \mu^{X_i}(1-\mu)^{1 - X_i}
      $$  is
 just the probability that
$$\Xi_{\delta,\epsilon}\left ((\bigcap_{i\in K(x)}(F^{-1}_i \bigcup_{j\in
m}[\frac{j}{m},\frac{j+X_j}{m}])\cup ([0,1]\setminus C_i))_x, E_x\right )$$ holds
and the choice of $m$ and the fact that $x\in Z$ guarantee that this
probability is greater than $1 - \rho^2/2$.
Hence $\alpha_m \geq \lambda(Z)(1 - \rho^2)$.

           Now let  $p $ be the probability that
$\Gamma(X_0,X_1,\ldots,X_m)
\geq (1 - \rho/2)\lambda(Z)$. Obviously, $p\lambda(Z) +
(1-p)(1-{\epsilon}/2)\lambda(Z) \geq \alpha_m \geq (1 - \rho^2)\lambda(Z)$.
Solving for $p$ yields that $p \geq 1 - 2{\rho} $. Since $m$ was chosen
so large that the probability that
$\lambda(\cup_{j\in m}[\frac{j}{m},\frac{j + X_j}{m}]) < 2^{-N}$ is greater than $
2{\rho}$,
there is at
least one $a \in W_N$ such that
 $\lambda(U) > (1 -
\rho/2)\lambda(Z)$ where $U$ is
the set of all $x\in Z$ such that
$$\Xi_{\delta,\epsilon(\delta,\mu,\eta,k)}\left ((\bigcap_{i\in
K(x)}(F^{-1}_i a)\cup ([0,1]\setminus C_i))_x, E_x\right )$$ holds.
Obviously $$\lambda(U\cup E_2) = \lambda(U) _ \lambda(E_2) \geq (1- \rho/2)\lambda(Z) +
\lambda(E_2)$$
$$\geq (1- \rho/2)(\lambda(E_1) - \frac{\rho\lambda(\pi_d(E))}{6(1 - \rho/2))}+ \lambda(E_2)$$
$$\geq  (1- \rho/2)(\lambda(E_1) + \lambda(E_2))   - (1- \rho/2)\frac{\rho\lambda(\pi_d(E))}{6(1 -
\rho/2)}\geq$$
$$ (1-\rho/2)(\lambda(\pi_d(E)) - 2 \frac{\rho\lambda(\pi_d(E))}{6(1 - \rho/2)}) -
       (1- \rho/2)\frac{\rho\lambda(\pi_d(E))}{6(1 -
\rho/2)} \geq (1-\rho)(\lambda(\pi_d(E)))$$
                                         as required.
\stopproof

\begin{theor} Suppose that
 $\{F_i\}_{i\in k}$ are continuous functions from $[0,1]^d$ to $[0,1]$, $\eta > 0$
$\delta > 0$, $N\in\omega$ and
 $\{A_i\}_{i\in k}$ are measurable subsets of $[0,1]^d$.
Then there is $\epsilon > 0$
\label{main} such that for each     closed subset
$E\subseteq [0,1]^d$, if $\Xi_{\delta,\eta}(E,E)$
holds
then
$$\Xi_{d\delta,\epsilon}(\bigcap_{i\in k}((A_i\cap
F_i^{-1}a)\cup ([0,1]^d\setminus A_i)),
 E)$$ also holds for some elementary set $a \in W_N$.
 \end{theor}
\proof  Proceed by induction on $d$ noting that if $d = 1$ then
 this follows directly from Lemma~\ref{finitelymanyfunctions}. So assume that the lemma has been
established for $d$ and that
$\{F_i\}_{i\in k}$ are continuous functions from $[0,1]^{d+1}$ to $[0,1]$, $\eta > 0$,
$\delta > 0$, $N\in\omega$ and
 $\{A_i\}_{i\in k}$ are measurable subsets of $[0,1]^{d+1}$.
Let $B_i = \{(x,y)\in [0,1]^d\times [0,1] : \lambda((F_i^{-1}\{y\})_x) > 0\}$ and note that
$$B_i = \{(x,y)\in [0,1]^d\times [0,1] :(\exists K \mbox{ compact})(\lambda(K) > 0 \AND K\subseteq
(F_i^{-1}\{y\})_x\}$$
and, because
$F$ is continuous,  the relation
$K\subseteq
(F_i^{-1}\{y\})_x$ is Borel. Moreover, so is the statement  $\lambda(K) > 0$ and so the set $B$ is
$\Sigma^1_1$ and hence, measurable. Let $B_i^* $ be the inverse image of $B_i$ under the mapping
$(x,y) \mapsto (x, F_i(x,y))$ or, in other words,
$(x,y) \in B_i^*$ if and only if $\lambda((F_i^{-1}\{F(x,y)\})_x) > 0$.  Since $B_i^*$ is clearly
measurable, it
follows that so is $C_i = [0,1]^{d+1} \setminus B_i^*$.

Now, for each $ i\ \in k$, let  $\{f_i^j
: j \in I_i\}$ enumerate a maximal collection of
functions such that
\begin{itemize}
\item $f_i^j : C_i^j \to [0,1]$ where $C_i^j \subseteq [0,1]^d$ is compact
\item  $f_i^j$ is continuous
\item $f_i^j \subseteq B_i$
\item if $x\in C_i^j\cap C_i^{j'}$ then $f_j(x) \neq f_{j'}(x)$
\item $\displaystyle \int_{C_i^j} \lambda((F_i^{-1}\{f_j(x)\})_x) dx > 0$.
\end{itemize}
The first thing to notice is that, for each $ i \in k$, such a family must be countable and
therefore, $I_i \leq \omega$ without loss of generality. To see this let
$E_i^j =\{(x,y)\in [0,1]^d\times [0,1]:
F_i(x,y) = f_i^j(x)\}$.
 If $j\neq {j'}$ then  $E_i^j \cap E_i^{j'} = \emptyset$ and, moreover,
$$\lambda(E_i^j) =  \int_{C_i^j} \lambda((F_i^{-1}\{f_j(x)\})_x) dx
> 0$$
for any $j \in I_i$. Hence the family of sets  $E_i^j$ is countable for each $ i \in k$.

Next, it must be shown that $$\sum_{j \in I_i}\int_{C_i^j} \lambda((F_i)^{-1}\{f_j(x)\})_x)dx =
\lambda(B_i^*)$$ so suppose not. Then it must be that
$\lambda(B_i^* \setminus \cup_{j\in I_i}E_i^j) > 0$. Since each
$f_i^j$ is continuous and $ I_i\leq \omega$, it follows that $B_i \setminus ( \cup_{j\in
I_i}f_i^j )$  is $\Sigma^1_1$.
Hence it is possible to use the von
Neumann selection theorem to find a function $f$
 such that the domain of $f$ is $\pi_d(B_i \setminus ( \cup_{j\in
I_i}f_i^j  )) $
and $f$ is measurable. Since $\pi_d(B_i \setminus ( \cup_{j\in
I_i}f_i^j  )) $
 is also equal to
$\pi_d(
B_i^* \setminus ( \cup_{j\in
d}E_i^j )) $  it must be that $\lambda(\pi_d(B_i \setminus ( \cup_{j\in
I_i}f_i^j  ))
) > 0$.
Hence $$\int_{\pi_d(B_i \setminus ( \cup_{j\in
I_i}f_i^j  )) } \lambda((F_i^{-1}\{f(x)\})_x) dx > 0$$ because
$\lambda(F_x^{-1}\{f(x)\}) > 0$ for each $ x\in \pi_d(B_i \setminus ( \cup_{j\in
I_i}f_i^j  )) $. Finally, by using Lusin's Theorem, it is
possible to find a compact set, $D\subseteq \pi_d(B_i \setminus ( \cup_{j\in
I_i}f_i^j  )) $ such that $f\restriction D$ is continuous and
$\int_{D} \lambda((F_i^{-1}\{f(x)\}_x) dx > 0$.
  This contradicts the putative maximality of the family $\{f_i^j : j \in I_i\}$.

Now note that for each $ i \in k$ the function $ (F_i \restriction
C_i)_x$ has small fibres for all $x$.
Applying Lemma~\ref{average} to $\{F_i\restriction (C_i\cap A_i) :i\in
k\}$, $\delta$,  $\eta/2$ and $N+1$
it follows that  there is some $\epsilon^* > 0$
such that for all $\mu > 0$ and any closed $E\subseteq [0,1]^{d+1}$ there is some  $a\in W_{N+1}$
such that the Lebesgue measure of
$$ \{x\in \pi_d(E) :
\Xi_{\delta,
\epsilon^*}((\bigcap_{i\in k}((C_i\cap A_i\cap F_i^{-1}a)\cup
([0,1]^{d+1}\setminus (C_i\cap A_i) ))_x, E_x)$$
$$\mbox{ or }
\neg \Xi_{\delta,
\eta}(E_x,E_x)mbox{ or }  H^r(E_x) \leq \delta
\}
$$ is at least $(1 - \mu)\lambda(E)$.

  It is therefore possible to find $K \in \omega$ such that for each $i\in k$
$$\sum_{j \in K}\int_{C_i^j} \lambda((F_i^{-1}\{f_i^j(x)\})_x)dx >
                                                                   \lambda(B^*_i) -
\frac{\eta^d\epsilon^*}{2^{d+1}k^2}$$ and so, if $S_i$ is defined to be
$$\{x\in \pi_d(E) : \lambda((B_i^* \setminus F_i^{-1}\{f_i^j(x)
: j\in K\})_x) \geq \epsilon^*/2k\}$$ then $\lambda(S_i) < \frac{\eta^d}{2^dk}$ for each $ i\in k$.
Let $U \subseteq [0,1]^d$ be any closed set  such that $U\cap S_i = \emptyset$ for each $ i\in
k$
and $\lambda(U) > 1 - (\eta/4)^d$.  Let $F_i^j$ be an arbitrary continuous extension of $f_i^j $
which has domain
$[0,1]^d$ and let $A_i^j = \dom(f_i^j)$.  It follows from the induction hypothesis
that
there is  $\epsilon' >0$
 such that if  $E\subseteq [0,1]^{d}$ is a closed set such that
$\Xi_{\delta,\eta/2}(E,E)$ holds then there is $a \in W_{N+1}$ such that
$$\Xi_{d\delta,\epsilon'}(\bigcap_{i\in k}\bigcap_{j\in K}((A_i^j\cap
(F_i^j)^{-1}a)\cup ([0,1]^d\setminus A_i^j)), E)$$ holds.
Let $\epsilon = \min\{\epsilon^*/2, \epsilon'/2, \eta/4\}$.

Now suppose that $E$ is a closed set such that $\Xi_{\delta,\eta}(E,E)$
holds. From the choice of  $\epsilon^*$ it follows that it is possible to find
  $a_0\in W_{N+1}$ such that
the Lebesgue measure of $Z=$
$$ \{x\in [0,1]^d :
\Xi_{\delta,
\epsilon^*}((\bigcap_{i\in k}((C_i\cap A_i\cap F_i^{-1}a_0)\cup
([0,1]^{d+1}\setminus (C_i\cap A_i) ))_x, E_x)$$
$$ \mbox{ or }\neg
\Xi_{\delta,
\eta}(E_x,E_x)   \mbox{ or } H^r(E_x) \leq \delta
\}
$$ is at least $1 - (\epsilon'/2)^d$.
If  $\hat{E} = \{x\in \pi_d(E) : \Xi_{\delta,\eta/2}(E_x,E_x) $ then $\Xi_{\delta,\eta}(\hat{E},
\hat{E})$ holds,  by Lemma~\ref{Xi.5}, because $\Xi_{\delta,\eta}({E},
{E})$ does.  From Lemma~\ref{Xi.2} it follows that $\hat{E}$ is Borel and so there exists a
closed set  $\bar{E}\subseteq \hat{E}$ such that $\lambda(\hat{E}\setminus \bar{E}) < (\eta/2)^d$.
Therefore  $\Xi_{\delta,\eta/2}({\bar{E}},
{\bar{E}})$ holds by lemma~\ref{ddimm}.
Another appeal to  Lemma~\ref{ddimm} yields that
$\Xi_{\delta,\eta/2}(\bar{E}\cap U,\bar{E})$ and so,  from Lemma~\ref{Xi.3} it may be
concluded  that $\Xi_{\delta,\eta/2}(\bar{E}\cap U,\bar{E}\cap U)$.

The choice of
$\epsilon'$ guarantees that
there is $a_1 \in W_{N+1}$ such that
$$\Xi_{d\delta,\epsilon'}(\bigcap_{i\in k}\bigcap_{j\in K}((A_i^j\cap
(F_i^j)^{-1}a_1)\cup ([0,1]^d\setminus A_i^j)), \bar{E}\cap U)$$ holds.
 It follows from Lemma~\ref{ddimm} that so does
$$\Xi_{d\delta,\epsilon'/2}(Z\cap \bigcap_{i\in k}\bigcap_{j\in K}((A_i^j\cap
(F_i^j)^{-1}a_1)\cup ([0,1]^d\setminus A_i^j)), \bar{E}\cap U)$$    and from Lemma~\ref{Xi.4} it
follows that
$$\Xi_{(d+1)\delta,\epsilon}(Z\cap \bar{E}\cap U\cap\bigcap_{i\in k}\bigcap_{j\in K}((A_i^j\cap
(F_i^j)^{-1}a_1)\cup ([0,1]^d\setminus A_i^j)), \pi_d(E))$$ because
$\Xi_{\delta, \epsilon}(\bar{E}\cap U,\pi_d(E))$ holds since $\epsilon \leq \eta/2$.
Let $a = a_0 \cup a_1$ and notice that $a \in
W_N$.

Using Lemma~\ref{reverse}, it suffices to show that if
$$x \in  Z\cap \bar{E}\cap U\cap \bigcap_{i\in k}\bigcap_{j\in K}((A_i^j\cap
(F_i^j)^{-1}a_1)\cup ([0,1]^d\setminus A_i^j))  $$
then
$$\Xi_{\delta,\epsilon}(\bigcap_{i\in k}((A_i\cap E\cap F_i^{-1}a)\cup ([0,1]^{d+1}\setminus
A_i))_x, E_x)$$ holds.
To see that this is so, recall that
since $x \in U$ it must be that
$\lambda(Y(x)) <  \epsilon^*/2$ where
$$Y(x) = \bigcup_{i\in k} (B_i^* \setminus F_i^{-1}\{f_i^j(x)
: j\in K\})_x) $$
Moreover, since $x \in Z$ it must be that
either
$$\Xi_{\delta,
\epsilon^*}(\bigcap_{i\in k}((C_i\cap A_i\cap F_i^{-1}a_0)\cup
([0,1]^{d+1}\setminus (C_i\cap A_i) ))_x, E_x)$$ holds or
$\Xi_{\delta,
\eta}(E_x,E_x)$ fails or $H^r(E_x) \leq \delta$.
However, since $x \in \bar{E}\subseteq \hat{E}$ it must be that $\Xi_{\delta,
\eta}(E_x,E_x)$  holds. If $H^r(E_x) \leq \delta$ then $\Xi_{\delta,
\eta}(\emptyset,E_x)$ holds and so, in either case it
                                                      follows
                                                             that
$$\Xi_{\delta,
\epsilon^*}((\bigcap_{i\in k}((C_i\cap A_i\cap F_i^{-1}a_0)\cup
([0,1]^{d+1}\setminus (C_i\cap A_i) ))_x, E_x)$$ holds.
          It therefore  follows from Lemma~\ref{ddimm} that
$$\Xi_{\delta,
\epsilon^*/2}((\bigcap_{i\in k}((\cap C_i\cap A_i\cap F_i^{-1}a_0)\cup
([0,1]^{d+1}\setminus (C_i\cap A_i))
))_x\setminus Y(x), E_x)\}$$ holds. Therefore it suffices to show that
           $$(E\cap C_i\cap A_i\cap F_i^{-1}a_0)\cup
(E\setminus (C_i\cap A_i))
))_x\setminus Y(x) \subseteq  ((A_i\cap E\cap F_i^{-1}a)\cup (E\setminus A_i))_x$$
for each $ i\in k$.

Fix $ i \in k$ and suppose that
$y \in (E\cap C_i\cap A_i\cap F_i^{-1}a_0)\cup
(E\setminus (C_i\cap A_i
))_x\setminus Y(x)$. If $y \in E\cap C_i\cap A_i\cap F_i^{-1}a_0$ then
$y  \in A_i\cap E\cap F_i^{-1}a$. On the other hand, suppose
$y \in (E\setminus (C_i\cap A_i))_x\setminus Y(x)$. If $y \in
E\setminus A_i$ there is nothing to prove so it may be assumed that
$y \in (A_i\setminus C_i)_x \setminus Y(x)$.
 Then,
since  $B^*_i\cap E = E\setminus C_i$ it must be that $y\in (B^*_i)_x$ and,
since    $y \notin Y(x)$, it follows that
$y\in F_i^{-1}\{f^j_i(x) : j\in K\}$ and so there is some $m \in K$
such that  $F_i(y) = f^m_i(x)$ and, in particular, $x \in A_i^m$.
 Since $$x \in \bigcap_{i\in k}
\bigcap_{j\in I_i}((A_i^j\cap
(F_i^j)^{-1}a_1)\cup ([0,1]^d\setminus A_i^j))$$
it follows that $x \in ((A_i^m\cap
(F_i^m)^{-1}a_1)$ and so $F_i(y) = F_i^m(x)\in a_1$.
Recalling that $y \in (E\setminus (C_i\setminus A_i))_x$  it follows
that
 $y \in (A_i \cap E\cap F_i^{-1}a)_x$.
\stopproof

\begin{corol}
For any $n\in \omega$ the ideal ${\cal I}^r_n$ is proper.
\end{corol}
\proof From Lemma~\ref{FiniteUnions} it suffices to show that that if
$X(f,{\cal C}, \delta')$ is a $d$-dimensional generator for ${\cal
I}^r_n$
then $W_n\not\subseteq X(f,{\cal C}, \delta')$. Let $\beta : \omega \to (0,1)$ witness that
${\cal C} = \{C_i\}_{i\in \omega}$ is a normal family. Let $ m$ be any integer such that $m >
d/\delta'$. Now apply Theorem~\ref{main} letting $\{F_i\}_{i\in k} = \{f\}$,
$\eta = \beta(m)$, $\delta = \delta'/d$,
$N= n$ and
$\{A_i\}_{i\in k} = \{[0,1]^d\}$, . This
yields $\epsilon >0$ such that for every $i\in \omega$ there is some $a_i \in
W_n$ such that
$\Xi_{d\delta,\epsilon}(f^{-1}a_i\cup([0,1]^d\setminus C_i, C_i)$ holds provided that
$\Xi_{\delta,\beta(m)}(C_i,C_i)$ holds. Without loss of generality, it may be assumed that
$\epsilon \leq \beta(m)$. This implies that
$\Xi_{d\delta,\epsilon}(f^{-1}a_i, C_i)$ holds provided that
$\Xi{\delta,\beta(m)}(C_i,C_i)$ does.
Since
$1/m < \delta$ and $\Xi_{1/m,\beta(m)}(C_i,C_i)$ holds for each $i \geq m$ it follows
$\Xi_{d\delta,\epsilon}(f^{-1}a_j, C_j)$ holds
for some $j$ such that  $\lambda(\pi_n(C_j) \setminus \pi_n(\cap_{i\in \omega}C_i)) <
(\frac{\epsilon}{d+1})^n$ for all $n \leq d$. The fact that such a $j$ exists follows from the
remark after Definition~\ref{normfam}.
Therefore
$\Xi_{\delta',\epsilon}(f^{-1}a, C_i)$ holds for all $i  > m$ by Corollary~\ref{equivs3}
and hence $a \notin X(f,{\cal C}, \delta')$.
\stopproof
\section{The End}
Finally, everything must be put together.
\begin{theor}
If ${\cal I} = \{{\cal I}^r_n\}_{n\in \omega}$ then
$$1\forces{\poset({\cal I })}{\lambda(V\cap [0,1]) = 0}$$ where $V$
represents the ground model.\end{theor}
\proof First notice that if $A \in {{\cal I}^r_n}^+$ then
$[0,1]\subseteq \cup X$  because if $x \in [0,1]$ then, letting
$\hat{x} : [0,1] \to [0,1]$ represent the function which is constantly $x$, it follows
that
$A\not\subseteq X(\hat{x}, \{[0,1]\}_{i\in\omega},1/2)$.  A standard genericity argument will
yield that if $G\in \prod_{n\in \omega}W_n$ is obtained from a $\poset({\cal I }) $ generic set
and $x \in V$ then  $x \in G(n)$ for infinitely any $n$.
Since any member of $W_n$ has measure less than $2^{-n}$, the result  is proved.
\stopproof
\makeatletter \renewcommand{\@biblabel}[1]{\hfill#1.}\makeatother
\renewcommand{\bysame}{\leavevmode\hbox to3em{\hrulefill}\,}

\end{document}